\documentclass[a4paper,10pt]{article}
\usepackage{geometry}                
\usepackage[active]{srcltx}
\usepackage{hyperref}
\usepackage[normalem]{ulem}

\usepackage{graphicx}
\usepackage{amssymb}
\usepackage{epstopdf}
\usepackage{amsmath}
\usepackage{amsthm}
\usepackage{amssymb}
\usepackage{latexsym}
\usepackage{mathtools}
\usepackage{mathrsfs}
\usepackage{graphics}
\usepackage{latexsym}
\usepackage{psfrag}
\usepackage{import}
\usepackage{verbatim}
\usepackage{enumerate}
\usepackage{enumitem}
\usepackage{graphicx}
\usepackage[usenames]{color}

\newcommand{\R}{\mathbb{R}}
\newcommand{\N}{\mathbb{N}}

\newcommand{\cH}{\mathcal{H}}

\newcommand{\LL}{\text{\rm L}}

\newcommand{\supp}{\text{\rm supp}}

\newcommand{\weak}{\rightharpoonup}

\newcommand{\T}{\mathcal{T}}

\newcommand{\CD}{\mathsf{CD}}

\newcommand{\Geo}{{\rm Geo}}
\newcommand{\TGeo}{{\rm TGeo}}
\newcommand{\MCP}{\mathsf{MCP}}

\newcommand{\Ent}{{\rm Ent}}

\newcommand{\Dom}{{\rm Dom}}
\newcommand{\fs}{\mathfrak{s}}
\newcommand{\fc}{\mathfrak{c}}
\newcommand{\fI}{\mathfrak{I}}
\newcommand{\cK}{\mathcal{K}}
\newcommand{\fa}{\mathfrak{a}}

\newcommand{\mui}{\mu_\infty}

\newcommand{\Prob}{\mathcal P}

\newcommand{\Ric}{{\rm Ric}}

\newcommand{\tr}{\textrm{tr}}
\newcommand{\K}{\mathcal{K}}

\newcommand{\vol}{\mathrm{vol}}

\newcommand{\TMCP}{\mathsf{TMCP}}
\newcommand{\TCD}{\mathsf{TCD}}
\newcommand{\wTCD}{\mathsf{wTCD}}

\renewcommand{\H}{\mathcal{H}}

\newcommand{\mm}{\mathfrak m}
\newcommand{\qq}{\mathfrak q}

\newcommand{\ee}{{\rm e}}

\newcommand{\sfd}{\mathsf d}

\renewcommand{\cH}{\mathcal{H}}

\theoremstyle{plain}
\newtheorem{lemma}{Lemma}[section]
\newtheorem{theorem}[lemma]{Theorem}

\newtheorem{proposition}[lemma]{Proposition}
\newtheorem{corollary}[lemma]{Corollary}

\newtheorem*{theorem*}{Theorem}
\newtheorem*{maintheorem*}{Main Theorem}

\theoremstyle{definition}

\newtheorem{definition}[lemma]{Definition}
\newtheorem*{definition*}{Definition}
\newtheorem{remark}[lemma]{Remark}

\numberwithin{equation}{section}

\title{A review of  Lorentzian synthetic theory \\of timelike Ricci curvature bounds}
\author{Fabio Cavalletti\thanks{F. Cavalletti: Mathematics Area, SISSA, Trieste (Italy), email: cavallet@sissa.it. } \ and Andrea Mondino\thanks{A. Mondino (corresponding author): Mathematical Institute,  University of Oxford  (UK), St. Hilda's College.  email: Andrea.Mondino@maths.oxford.ac.uk.}}

\date{\today}     

\begin{document}
\maketitle

\begin{abstract}
The goal of this survey is to give a self-contained introduction 
to synthetic timelike Ricci 
curvature bounds for (possibly non-smooth) Lorentzian spaces via 
optimal transport \& entropy tools, including a synthetic version of Hawking's 
singularity theorem and a synthetic characterisation of Einstein's vacuum equations.\ 
We will also discuss some motivations arising from the smooth world
and some possible directions for future research.
\end{abstract}

\bibliographystyle{plain}

\tableofcontents

\section*{Introduction}

Optimal transport  turned out to be a very effective tool to study spaces with extremely low regularity in Riemannian signature.\ 
More precisely, one can use optimal transport to define and study metric measure spaces with Ricci $\geq K\in \R$ and dimension $\leq N\in[1,\infty)$ in a synthetic sense; 
this is indeed the celebrated theory of $\CD(K,N)$ spaces pioneered by Sturm \cite{sturm:I, sturm:II} and Lott-Villani \cite{lottvillani:metric} independently.

Motivated by the success of optimal transport techniques in the Riemannian signature, during the last few years the development of optimal transport tools in the Lorentzian setting has become an increasingly popular topic of research, let us give a brief account.

The optimal transport problem for relativistic costs was proposed by Brenier \cite{Br03}, further investigated in \cite{BertrandPuel, BertrandPratelliPuel}, and pushed to a geometric Lorentzian setting in \cite{EM17,Suhr, KellSuhr}.\  An intriguing physical motivation for studying optimal transport in a Lorentzian setting is the \emph{early universe reconstruction problem}.  Studied in \cite{BFHLMMS03} and in the Nature paper \cite{FMMS02}  with methods
of optimal transportation, such a problem amounts to reconstructing the trajectories of masses
from the big bang to their present day positions in Robertson-Walker spacetimes.\
A mathematical formulation for general globally hyperbolic spacetimes would read
as follows: given two probability  measures, one concentrated on a Cauchy hypersurface, the
other on the past cone of a point, what can be said about the trajectories
of the minimizers in a dynamical optimal coupling of the two measure? In \cite{FMMS02} it is explained why the problem can be attacked with methods of optimal
transportation.
\\

Due to the central role played by the Ricci curvature in general relativity (e.g.   Einstein's equations can be formulated just in terms of the Ricci tensor, the energy momentum tensor, and the cosmological constant) and due to the success of optimal transport tools to study Ricci curvature bounds in the Riemannian signature, there has been a growing interest in the interplay of optimal transport and Ricci curvature bounds in a Lorentzian setting.\ Indeed, for a smooth Lorentzian manifold (see Section \ref{sec:TRLBsmooth}),  McCann \cite{McCann} characterised non-negative timelike Ricci curvature lower bounds in terms of dispacement convexity, and the second author joint with Suhr \cite{MoSu} characterised any lower and upper bound on the timelike Ricci curvature (and thus the full Einstein equations) in terms of optimal transport and entropy.
\\

While the non-smooth (or “synthetic") counterpart of a Riemannian manifold is a metric (measure) space, and thus the latter gives a natural framework for Lott-Sturm-Villani theory, the interplay of causal and metric structures in a Lorentzian manifold is more subtle and forces to modify the setting in order to develop a synthetic theory.\  Such a program was pioneered in seminal work by Kronheimer-Penrose \cite{KroPen67} and has been recently  formalised  by Kunzinger-S\"amann \cite{KS} into the concept of Lorentzian (pre-)length spaces (see Section \ref{sec:LSS}).\ 
Lorentzian (pre)-length spaces give a natural framework to develop a theory of optimal transport and of synthetic curvature bounds.
\\

Synthetic \emph{timelike sectional curvature bounds} for Lorentzian (pre)-length spaces have been investigated in \cite{KS} (see also the more recent \cite{BeSa}).\ Synthetic timelike Ricci curvature lower bounds in Lorentzian (pre)-lenght spaces, which constitute the subject of the present survey, have been developed by the authors in \cite{CaMo:20} (inspired by the aforementioned smooth characterisations obtained in \cite{McCann, MoSu}).\ 
The aim of \cite{CaMo:20} has been to initiate 
a theory of Lorentzian $\CD$ spaces, denoted by $\TCD$ for ``timelike $\CD$'' (see Section \ref{Sec:Ricci>Synt}), and to 
obtain some first 
geometric and functional insights; for instance the following inequalities are among the consequences of the $\TCD$ condition: volume comparison, diameter bound and  Poincar\'e inequality.\ 
Moreover, a far reaching synthetic version of the Hawking's singularity theorem holds (in  a sharp form, and in a formulation which includes Lorentzian metrics of very low regularity: $C^{0,1}$ or, more generally, $C^{0}$ plus causally plain).\ We refer to Section \ref{S:applications} for an overview of (some of) the applications.
\\

In Section \ref{S:SEE} we survey synthetic timelike Ricci upper bounds and a synthetic formulation of Einstein's vacuum equations for Lorentzian (pre)-length spaces, obtained in \cite[App. B]{MoSu} after \cite{StUB, CaMo:20}.\
We will end the review  by proposing some possible directions for future research, see Section \ref{Sec:Open}.
\\

Let us conclude the introduction  by stressing that, apart from the specific synthetic framework and the aforementioned results, we expect optimal transport tools to be useful in a wide range of applications.\ Indeed, typically, classical arguments make use of Jacobi fields computations (e.g. Raychaudhuri equation) heavily relying on the $C^{2}$-regularity of the metric (with some effort, one can lower to $C^{1,1}$ by approximation).\
The main advantage of the optimal transport tools surveyed here (see \cite{CaMo:20} for more results and details) is that they often allow to carry over arguments in a setting of very low regularity (including $C^{0,1}$ or, more generally, $C^{0}$ and causally plain) where the classical terms are even not well defined.
\\

 Motivations for studying Lorentzian metrics of low regularity (or lowering the regularity even to a Lorentzian synthetic framework) come both from the PDE point of view in general relativity (i.e. the Cauchy initial value problem for the Einstein equations) and from physically relevant models. 

From the PDE point of view, the standard local existence results for the vacuum Einstein equations assume the metric to be of Sobolev regularity $H^{s}_{loc}$, with $s>\frac{5}{2}$ (see for instance \cite{Rendall}).\ The Sobolev regularity of the metric has been lowered even further (e.g. \cite{KRS}).\ Related to the initial value problem for the Einstein equations, one of the main open problems in the field is the so called (weak/strong) censorship conjecture proposed originally by Penrose and later refined in \cite{Chri} (see also \cite{Daf} for a more updated overview).\ The strong form of the conjecture asserts (roughly) that the maximal globally hyperbolic development of generic initial data for the Einstein equations is inextendible as a suitably regular Lorentzian manifold.\
Formulating a precise statement of the conjecture is itself non-trivial since one needs to give a precise meaning to ``generic initial data'' and ``suitably regular Lorentzian manifold''.\ Understanding the latter is where Lorentzian metrics of low regularity and related inextendibility results play a pivotal role.\  The strongest form of the conjecture would prove inextendibility for a  $C^{0}$ metric.\  As pointed out by Chrusciel-Grant \cite{CG} (see also \cite{Ling} for a survey), causality theory for $C^{0}$ metrics departs significantly from classical theory (e.g. the lightlike curves emanating from a point may span a set with non-empty interior, a phenomenon called ``bubbling'').\ Nevertheless, Sbierski \cite{Sbi} gave a clever proof of $C^{0}$-inextendibility of Schwarzschild, Minguzzi-Suhr  \cite{MiSu} showed $C^{0}$-inextendibility for timelike geodesically complete spacetimes, and Grant-Kunzinger-S\"amann \cite{GKS} pushed the inextendibility to Lorentzian length spaces.

As already envisaged in seminal work by Lichnerowitz in the 50'ies \cite{Lichn}, from the point of view of physically relevant models, several types of matter in a spacetime may give a discontinuous energy-momentum tensor and thus, via the Einstein equations, lead to a Lorentzian metric of regularity lower than $C^{2}$.\ Examples of such a phenomenon are spacetimes that model the inside and outside of a star,  matched spacetimes \cite{MaSe}, self-gravitating compressible fluids \cite{BuLe} or shock waves.\ Some physically relevant models require even lower regularity, for instance: spacetimes with conical singularities \cite{Vick2},  cosmic strings \cite{Vick1} and (impulsive) gravitational waves (see for instance the pioneering work of Penrose \cite{PenGW}, and the more recent \cite[Chapt. 20]{GrPo}), or models for cyclic spacetimes \cite{LLV}.

Finally, a long term motivation for studying non-regular Lorentzian spaces is the desire of understanding the ultimate nature of spacetime.\ The rough picture is that at the quantum level (and thus in extreme physical conditions such as gravitational collapse or the origin of the universe), the spacetime may be very singular and possibly not approximable by smooth  structures (see the end of Section \ref{sec:LSS}). 

In case of a metric of low regularity, the approach to curvature used so far is distributional, exploiting the smoothness of the underline spacetime.\ This allows  \cite{GeTr} (see also \cite{StVick}) to define distributional curvature tensors for  $W^{1,2}_{loc}$-Lorentzian metrics  satisfying a suitable non-degeneracy condition (fulfilled for instance when the metric is $C^{1}$, see  \cite{Graf}).\ 
One of the goals of the approach surveyed here  is to address the question of (timelike Ricci) curvature when \emph{not only the the metric tensor, but the space itself is singular}.

\subsubsection*{Acknowledgements}
The second  author is  supported by the European Research Council (ERC), under the European's Union Horizon 2020 research and innovation programme, via the ERC Starting Grant CURVATURE, grant agreement No. 802689.

\section{Optimal Transport in Lorentzian spaces}\label{S:Basics}

\subsection{Classical setting and Lorentzian length spaces} \label{sec:LSS}
We start by introducing some basic terminology. 

Given  a Lorentzian manifold $(M^{n},g)$, we also fix an auxiliary complete Riemannian metric $h$ on $M$ and we denote by $\sfd^{h}$ the associated distance.\ 
We say that  $v\in T_{x}M\setminus \{0\}$ is \emph{timelike} if $g(v,v)< 0$.\ 
$(M,g)$ is \emph{time-oriented} if it admits a continuous timelike vector field $Y$.\ Such a
$Y$ induces a partition on the set  of 
timelike vectors: the equivalence class of 
\emph{future pointing} vectors  $\{v\in  T_{x}M  \mid  g(v,Y) < 0\}$ and the class of \emph{past pointing} vectors 
 $\{v\in  T_{x}M  \mid  g(v,Y) > 0\}$.\ 
 
 The closure of the set of future pointing timelike vectors is called \emph{future causal cone} and denoted by $\mathcal{C}\subset TM$.\ Let $I=[0,1]$ denote the unit interval.\
A locally Lipshitz curve $\gamma:I\to M$ 
is  \emph{causal} if $\dot \gamma \in \mathcal{C}$ for a.e. $t$; 
\emph{timelike} if $\dot \gamma \in  {\rm Int}(\mathcal{C})$ for a.e. $t$.
\\These definitions induce two causal relations 
for $x,y \in  M$: 
\begin{itemize}
\item we say that 
$x \leq y$ if $\exists \ \gamma$  causal curve with 
$\gamma_0 = x$ and $\gamma_1 = y$; 
\item we say that  $x \ll y$ 
if $\exists \ \gamma$  timelike curve with 
$\gamma_0 = x$ and $\gamma_1 = y$. 
\end{itemize}
The \emph{length} of a causal curve $\gamma: I\to M$ is defined as
$$
{\rm L}_{g}(\gamma):= \int_{I} \sqrt{|g(\dot \gamma,\dot \gamma)|}\, {\rm d}t. 
$$
The \emph{time separation} function $\tau : M \times M \to [0, \infty]$ is defined as
\begin{equation}\label{def:tausmooth}
\tau(x,y): = 
\begin{cases}
\sup \, \{ {\rm L}_{g}(\gamma): \gamma \text{ causal,  } \gamma_{0}=x, \, \gamma_{1}=y\}, & \text{ if } x \leq y, \\
0,& {\rm otherwise}.
\end{cases}
\end{equation}
The following reverse triangle inequality holds true for $\tau$: 
\begin{equation}\label{eq:RTI}
\tau (x , y ) + \tau (y , z ) \leq \tau (x , z ),\quad \text{ if } x \leq y \leq z. 
\end{equation}
If $g$ is $C^{0,1}$-regular ($C^{0}$ is not enough, but $C^{0}$ and causally plain suffices \cite{KS} after \cite{CG}):
$\tau$ is lower-semicontinuous and $\tau(x,y) > 0$ if and only if
 $x \ll y$.
 
The idea to formulate a synthetic theory for Lorentzian spaces is to turn the previous properties into definitions.\ Such a point of view has roots in the seminal work of Kronheimer-Penrose \cite{KroPen67} and has been formalised in recent work by Kunzinger-S\"amann \cite{KS}.

\begin{definition}[Causal space  $(X,\ll,\leq)$]
A \emph{causal space}  $(X,\ll,\leq)$ is a set $X$ endowed with a preorder $\leq$ and a transitive relation $\ll$ contained in $\leq$.
\end{definition}

One says that $x$ and $y$ are \emph{timelike} (resp. \emph{causally}) related if $x\ll y$  (resp. $x\leq y$).\  Let $A\subset X$ be an arbitrary subset of $X$.\ We define the \emph{chronological} (resp. \emph{causal}) future of $A$ the set
$$
I^{+}(A):=\{y\in X\mid \exists x\in A,\, x\ll y\}, \quad 
J^{+}(A) :=\{y\in X\mid \exists x\in A,\, x\leq y\}
$$
respectively.\ Analogously, one can define the   \emph{chronological} (resp. \emph{causal}) past of $A$. 
In order to keep the notation short, it is also useful to set
$$
X^2_{\leq}:=\{(x,y)\in X^2\mid x\leq y\}, \quad X^2_{\ll}:=\{(x,y)\in X^2\mid x\ll y\}.
$$

Additionally requiring the existence of a time separation function 
produces the following  notion that has to be understood like 
the Lorentzian analog of metric spaces.

\begin{definition}[Lorentzian pre-length space $(X,\sfd, \ll, \leq, \tau)$]
A \emph{Lorentzian pre-length space} \\ $(X,\sfd, \ll, \leq, \tau)$ is a  casual space $(X,\ll,\leq)$ additionally  equipped with a proper metric $\sfd$  (i.e. closed and bounded subsets are compact) and a lower semicontinuous function $\tau: X\times X\to [0,\infty]$,  called \emph{time-separation function}, satisfying
\begin{equation}\label{eq:deftau}
\begin{split}
\tau(x,y)+\tau(y,z)\leq \tau (x,z) &\quad\forall x\leq y\leq z \quad \text{reverse triangle inequality} \\
\tau(x,y)=0, \; \text{if } x\not\leq y, & \quad  \tau(x,y)>0 \Leftrightarrow x\ll y.
\end{split}
\end{equation}
\end{definition}

Notice that $X$ is endowed with the metric topology induced by 
$\sfd$ and all the topological concepts on $X$ (like l.s.c. of $\tau$) are formulated in terms of this topology.
Note that the lower semicontinuity of $\tau$ implies that $I^{\pm}(x)$ is open, for any $x\in X$.

Analogously to the causal relations, 
$\gamma: I \to X$ is said \emph{timelike} (resp. \emph{causal}) if it is locally Lipschitz (w.r.t. $\sfd$) and $\gamma_{t_{1}} \ll \gamma_{t_{2}}$ (resp. $\gamma_{t_{1}} \leq \gamma_{t_{2}}$) for all $t_{1}<t_{2}$.\
The length induced by $\tau$ for a causal curve $\gamma$  is then defined as: 
$$
{\rm L}_{\tau}(\gamma): = \inf \,  \sum_{i} \tau(\gamma_{t_{i}},\gamma_{t_{i+1}}) ,
$$
where the $\inf$ is taken over all finite partitions $0=t_{0}<t_{1}<\ldots<t_{m}=1$, $m\in \N$, of $I=[0,1]$.
\\A causal curve $\gamma : [a,b] \to X$ is a \emph{geodesic} if 
${\rm L}_{\tau}(\gamma) = \tau(\gamma_{a},\gamma_{b})$.

\noindent
We denote the set of causal (resp. timelike) geodesics as:
\begin{align}
&\Geo(X):=\{ \gamma:[0,1]\to X\,:  \, \tau(\gamma_{s}, \gamma_{t})=(t-s) \tau(\gamma_{0}, \gamma_{1})\, \forall s<t\}, \\
& \TGeo(X):=\{ \gamma \in \Geo(X):  \, \tau(\gamma_{0}, \gamma_{1})>0\}. \label{eq:defTGeo}
\end{align}
Note that, by definition, geodesics are always \emph{maximising}, \emph{future oriented,} and \emph{parametrised at constant speed on $[0,1]$}.

Moreover a Lorentzian pre-length space $(X,\sfd,\ll,\leq,\tau)$ is called:
\begin{itemize}
\item \emph{Causally closed:} $\{x\leq y \} \subset X\times X$ is closed
\item \emph{Non-totally imprisoning} if $\forall\, K \subset X$ compact, 
$\exists \, C >0$ s.t. ${\rm L}_{\sfd}(\gamma) \leq C$ for all $\gamma$ causal curve in $K$.
\item \emph {Globally hyperbolic} if it is non-totally imprisoning and 
$\forall\,x,y \in X$, $J^{+}(x) \cap J^{-}(y)$ is compact.

\item \emph{$\mathcal{K}$-Globally hyperbolic:}
if it is non-totally imprisoning and 
$\forall\,K_{1},K_{2} \subset X$ compact, $J^{+}(K_{1}) \cap J^{-}(K_{2})$ is compact.
\item \emph{Geodesic}: if $\forall\,x,y \in X$ with $x\leq y$ there exists   a  geodesic $\gamma$ from $x$ to $y$.
\end{itemize}

It was proved in \cite[Theorem 3.28]{KS} that for a globally hyperbolic  Lorentzian geodesic (actually length would suffice) space $(X,\sfd, \ll, \leq,\tau)$, the time-separation function $\tau$ is finite and continuous.\ Moreover, any globally hyperbolic  Lorentzian length space (for the definition of Lorenzian length space see \cite[Definition 3.22]{KS}, we skip it for brevity since we will not use it) is geodesic  \cite[Theorem 3.30]{KS}.

\subsection*{Examples entering the class of Lorentzian synthetic  spaces}

\noindent 
\textit{Spacetimes with a continuous Lorentzian metric.\ } Let $M$ be a smooth manifold endowed with a continuous Lorentzian metric $g$.\ Assume that $(M,g)$ is time-oriented.\  Observe that, for $C^0$-metrics, the natural class of differentiability of the underlying manifolds is $C^{1}$; now, $C^{1}$ manifolds always
admit a $C^{\infty}$ subatlas, and one can pick such a sub-atlas whenever convenient.\  We endow such a spacetime $(M,g)$ with the time separation function $\tau$ defined in \eqref{def:tausmooth}.

For a  spacetime with a Lorentzian $C^{0}$-metric: 
\begin{itemize}
\item Global hyperbolicity implies causal closedness and  $\cK$-global hyperbolicity  \cite[Prop. 3.3 and Cor. 3.4]{SaC0}.
\item Recall that a Cauchy hypersurface is a subset  which intersects exactly once every inextendible causal curve.\ It was proved in \cite{SaC0} that every Cauchy hypersurface is a closed acasual topological hypersurface and that global hyperbolicity is equivalent to the existence of a Cauchy hypersurface.
\item  If $g$ is a causally plain (or, more strongly, locally Lipschitz) Lorentzian $C^{0}$-metric on $M$ then the associated synthetic structure is a  Lorentzian pre-length space, see \cite[Prop. 5.8]{KS}.
More strongly, if $g$ is a   globally hyperbolic and causally plain Lorentzian $C^{0}$-metric on $M$ then the associated synthetic structure is a  causally closed,  $\cK$-globally hyperbolic Lorentzian geodesic space (see  \cite[Thm. 3.30 and Thm. 5.12]{KS}).
\end{itemize}
Summarising, if $M$ is a smooth manifold endowed with a $C^{0}$-Lorentzian metric $g$ making $(M,g)$ time-oriented then $(M,\ll,\leq)$ is a causal space 
but it is not necessarily a Lorentzian pre-length space. 
If $g$ is ``causally plain'' \cite{CG} (for instance 
$g\in C^{0,1}$), then $(M,\sfd^{h},\ll,\leq,\tau)$ is a Lorentzian pre-length space.\
If in addition $g$ is globally hyperbolic, then $(M,\sfd^{h},\ll,\leq,\tau)$ is a
causally closed, $\mathcal{K}$-globally hyperbolic and  geodesic Lorentzian space. 

Thus, the framework of causally closed, $\mathcal{K}$-globally hyperbolic and geodesic Lorentzian spaces is rather natural and will provide the setting of our work.
\medskip

\noindent 
\textit{Closed cone structures}.\  Closed cone structures can be seen as the synthetic-Lorentzian analogue of Finsler manifolds.\ They provide
a rich source of examples of Lorentzian pre-length and length spaces (see \cite[Sec. 5.2]{KS} for more details).\ We refer to  Minguzzi's comprehensive paper \cite{Min} for a thorough analysis of causality theory in the framework of closed cone structures, including embedding and singularity theorems.
\medskip

\noindent \textit{Some examples towards  quantum gravity}.\ The general setting of Lorentzian pre-length spaces allows to consider more general structures than Lorentz(-Finsler) metrics on smooth manifolds.  A remarkable motivation for such a general framework is given by certain approaches to quantum gravity.\ For instance let us mention  \cite{MP} where it is shown that, in a purely order theoretic manner,  one can reconstruct a globally hyperbolic spacetime and the causality relation from  a countable dense set of events.   
\\Two approaches to quantum gravity,  particularly close in spirit to Lorentzian pre-length spaces,  are the theory of \emph{causal Fermion systems} \cite{FinsterPrimer, FinsterBook} and the  \emph{theory of causal sets} \cite{BLMS}.\ The basic idea in both the approaches is that the structure of spacetime needs to be relaxed
on a microscopic scale to include quantum effects.\ This leads to non-smoothness of
the underlying geometry, and  the classical  structure of a Lorentzian spacetime emerges only in the macroscopic regime.\   For the connection to the theory of Lorentzian (pre-)length spaces, the reader is referred to  \cite[Sec. 5.3]{KS} and \cite[Sec. 5.1]{FinsterPrimer}.

\subsection{Optimal transport in Lorentzian spaces}

The space of Borel probability measures (resp. with compact support) over a metric space $(X,\sfd)$ is denoted by $\mathcal{P}(X)$ (resp. $\mathcal{P}_c(X)$).\ 
For $\mu,\nu\in \mathcal{P}(X)$, we denote by $\mu\otimes \nu\in \mathcal{P}(X^2)$ the unique probability measure such that
\begin{equation}\label{eq:defmutimesnu}
\mu\otimes \nu (A\times B):=\mu(A)\cdot \nu(B), \quad \text{for all $A,B\subset X$ Borel subsets.}
\end{equation}
Note that \eqref{eq:defmutimesnu} uniquely defines a probability measure on $X^2$, as the $\sigma$-algebra of Borel sets of $X^2$ is generated by products of Borel sets of $X$.

Let $P_{i}:X\times X\to X$, $i=1,2$, denote the projection maps on the factors. 
For a Borel map $f:X\to Y$, the associated push-forward map $f_{\sharp}:\mathcal{P}(X)\to \mathcal{P}(Y)$ is defined by 
$$
f_{\sharp}\mu (B):=\mu (f^{-1}(B)), \quad \text{for all $\mu\in \mathcal{P}(X)$ and $B\subset Y$ Borel subset.} 
$$
Notice that $(P_1)_{\sharp}(\mu\otimes \nu)=\mu$ and $(P_2)_{\sharp}(\mu\otimes \nu)=\nu$.

\begin{definition}[Transport plans]
If $(X,\sfd,\ll,\leq,\tau)$ is a  Lorentzian pre-length space 
 and $\mu,\nu \in \mathcal{P}(X)$, we can distinguish different families of
transport plans:
\begin{itemize}
\item[-]Classic: $\Pi(\mu,\nu) : = \{\pi \in \mathcal{P}(X^{2})  \mid (P_{1})_{\sharp}\pi = \mu, (P_{2})_{\sharp}\pi = \nu\}$. 
\item[-] Causal: $\Pi_{\leq}(\mu,\nu) : = \{ \pi \in \Pi(\mu,\nu) \mid \pi(X^{2}_{\leq}) = 1 \}$, %
\item[-] Timelike: $\Pi_{\ll}(\mu,\nu) : = \{ \pi \in \Pi(\mu,\nu) \mid \pi(X^{2}_{\ll}) = 1 \}$.
 \end{itemize}
\end{definition}

\textbf{Some intuition behind the notion of transport plan}.\
Since the survey is meant for an audience possibly not specialized in optimal transport, let us discuss some basic heuristics.
 
The rough idea of a classical transport plan $\pi\in \Pi(\mu,\nu)$ is that the ``mass $\mu({\rm d}x)$ at $x$'' is possibly split and transported following the law $\pi({\rm d}x {\rm d}{y})$.\   Note that if $f:X\to X$ is a Borel map, then $\pi:=({\rm Id}, f)_{\sharp}\mu\in \Pi(\mu, f_{\sharp}\mu)$ is a classical transport plan from $\mu$ to  $f_{\sharp}\mu$.\ However this is a very special case of a transport plan, namely a transport plan induced by a map (called transport map).\ A (trivial) example of a transport plan \emph{not} induced by a map is $\pi:=\mu\otimes\nu$ obtained by taking the product of $\mu$ and $\nu$ as in \eqref{eq:defmutimesnu}.\ It is also instructive to consider the case of $\mu=\delta_{x_{0}}$, $\nu=\frac{1}{2}\left( \delta_{x_{1}}+\delta_{x_{2}}\right)$, where $\delta_{x}$ denotes the Dirac mass at $x\in X$: in this case it is not possible to find a trasport map, and a transport plan is given by the product $\mu\otimes\nu$, which corresponds to ``splitting'' the Dirac mass at $x_{0}$ in half and transporting each half into $x_{1}$ and $x_{2}$ respectively.  

The rough idea of a causal (resp. timelike) transport plan $\pi\in \Pi_{\leq}(\mu,\nu)$ is that the ``mass $\mu({\rm d}x)$ at $x$'' is possibly split and transported following the law $\pi({\rm d}x {\rm d}{y})$ so that ``the destination is in the causal (resp. timelike) future of the source''.\ Clearly, for $\mu=\delta_{x_{0}}$ and $\nu=\delta_{x_{1}}$, it holds that $\pi=\mu\otimes\nu\in \Pi_{\leq}(\mu,\nu)$ (resp. $ \Pi_{\ll}(\mu,\nu)$) if and only if $x_{0}\leq x_{1}$ (resp. $x_{0}\ll x_{1}$).
\hfill$\Box$

\begin{definition}\label{def:Wp}
Let  $(X,\sfd, \ll, \leq, \tau)$ be a Lorentzian pre-length space and let $p\in (0,1]$.\ Given $\mu,\nu\in \mathcal{P}(X)$, the \emph{$p$-Lorentz-Wasserstein distance} is defined by
\begin{equation}\label{eq:defWp}
\ell_{p}(\mu,\nu):= \sup_{\pi \in \Pi_{\leq}(\mu,\nu)} \left(  \int_{X\times X}  \tau(x,y)^{p} \, \pi({\rm d}x{\rm d}y)\right)^{1/p}.
\end{equation}
When $\Pi_{\leq}(\mu,\nu)=\emptyset$, we set $\ell_{p}(\mu,\nu):=-\infty$.
\end{definition}
Definition \ref{def:Wp} extends to Lorentzian pre-length spaces the corresponding
notion given in the smooth Lorentzian setting in \cite{EM17}  
(see also \cite{McCann, MoSu}, and \cite{Suhr} for $p=1$). 
When $\Pi_{\leq}(\mu,\nu)=\emptyset$, we adopt the convention of McCann  \cite{McCann} (note that \cite{EM17} set $\ell_{p}(\mu,\nu)=0$ in this case).

The function $\ell_{p}$ inherits the properties of the time-separation function 
$\tau$.\ Indeed, by using the classical technique of gluing (used for instance to prove the triangle inequality for the Wasserstein distance), 
one can prove that $\ell_{p}$ verifies the reverse triangle inequality: 
\begin{equation}\label{E:triangular}
\ell_{p}(\mu_{0},\mu_{1}) +\ell_{p}(\mu_{1},\mu_{2}) \leq \ell_{p}(\mu_{0},\mu_{2}),\quad \forall\ \mu_{0},\mu_{1},\mu_{2}\in \mathcal{P}(X)
\end{equation}
with the convention on the left hand side that $\infty-\infty=-\infty$.

To invoke the classical theory of Optimal Transport it is more convenient to 
move the causal constraint on the transport plans to the cost function.\ In other words, it is useful to reformulate the variational problem on the right hand side of \eqref{eq:defWp}
with the following cost function: 
\begin{equation}\label{eq:defell}
\ell(x,y) : = 
\begin{cases}
\tau(x,y), & \textrm{if } x\leq y, \\
 -\infty, & \textrm{otherwise}.
\end{cases}
\end{equation}
Notice that
$$
 \int_{X\times X}  \tau(x,y)^{p} \, \pi({\rm d}x{\rm d}y)=   \int_{X\times X}  \ell(x,y)^{p} \, \pi({\rm d}x{\rm d}y)\in \R_{\geq 0}, \quad \text{ for all }\pi\in \Pi_{\leq}(\mu,\nu).
$$
Moreover, if $\pi\in \Pi(\mu,\nu)$ satisfies $ \int_{X\times X}  \ell(x,y)^{p} \, \pi({\rm d}x{\rm d}y)>-\infty$ then $\pi\in \Pi_{\leq}(\mu,\nu)$.\
Thus the maximization problem  \eqref{eq:defWp} is equivalent (i.e. the $\sup$ and the set of maximisers coincide) to the  maximisation problem
\begin{equation}\label{eq:supell}
\sup_{\pi \in \Pi(\mu,\nu)} \left(  \int_{X\times X}  \ell(x,y)^{p} \, \pi({\rm d}x{\rm d}y)\right)^{1/p}.
\end{equation}

A  $\pi\in  \Pi_{\leq}(\mu,\nu)$ maximising in \eqref{eq:defWp} is said to be \emph{$\ell_{p}$-optimal}.\ The set of \emph{$\ell_{p}$-optimal} plans from $\mu$ to $\nu$ is denoted by $  \Pi_{\leq}^{p\text{-opt}}(\mu,\nu)$.

The advantage of the formulation \eqref{eq:supell} is that, when $(X,\sfd, \ll, \leq, \tau)$ is causally closed (so that $\{(x,y) \mid x\leq y\}\subset X\times X$ is a closed subset) and globally hyperbolic geodesic (so that $\tau$ is continuous) then  $\ell$ is upper semi-continuous on $X\times X$.\
Similarly, when $X$ is  \emph{locally} causally closed globally hyperbolic geodesic, if $\mu$ and $\nu$ have compact support then  $\ell$ is upper semi-continuous on $\supp \, \mu \times \supp \, \nu$. 
\\In both cases, one can apply standard optimal transport techniques (e.g. \cite{Vil}) to the Monge-Kantorovich problem \eqref{eq:supell}.

In the following, given two functions $a,b:X\to \R$, we denote by $a\oplus b: X\times X\to \R$ the function on the product defined by
$$
(a\oplus b)(x,y):=a(x)+b(y), \quad \text{for all } x,y\in X.
$$

\begin{proposition}[Existence of optimal plans]\label{prop:ExMaxellp}
Let  $(X,\sfd, \ll, \leq, \tau)$ be a  causally closed  
(resp. locally causally closed) globally hyperbolic Lorentzian 
geodesic space and let $\mu,\nu\in \mathcal{P}(X)$ (resp.  $\mathcal{P}_{c}(X)$). 
If  $ \Pi_{\leq}(\mu,\nu)\neq \emptyset$ and if  there exist measurable functions $a,b:X\to \R$, with $a\oplus b \in L^{1}(\mu\otimes \nu)$ 
such that $\ell^{p}\leq a\oplus b$  on $\supp \, \mu \times \supp \, \nu$ 
(e.g. when $\mu$ and $\nu$ are compactly supported) then the $\sup$ in \eqref{eq:defWp} 
(and henceforth in \eqref{eq:supell}) is attained and finite. 
\end{proposition}

We fix the following notation 
\begin{align*}
\Pi^{p\text{-opt}}_{\leq} (\mu, \nu)&:=\{\pi \in \Pi_{\leq}(\mu,\nu) \text{ is  $p$-optimal}\},\\
\Pi^{p\text{-opt}}_{\ll}(\mu, \nu)&:=\{\pi \in \Pi_{\ll}(\mu,\nu)\text{ is  $p$-optimal}\}.    
\end{align*}

\subsection{The Kantorovich duality}

Existence of optimal plans $\pi \in \Pi_{\leq}(\mu,\nu)$ easily follows by 
the direct method of  the Calculus of Variations.\
The rough picture behind a transport plan is that  $\pi$ is moving the mass 
$\mu({\rm d}x)$ “at $x$” in other “suitable points” of $\supp\, \nu$.\
A natural question is then:  which are such “suitable points” chosen by $x$?
Or, better said, how can we detect $\Gamma \subset X^2_{\leq}$
such that $\pi(\Gamma) = 1$?

These questions can be answered via two key notions in optimal transport theory: 
cyclical monotonicity and Kantorovich duality.\
Optimal transport problems possess indeed a rich duality theory, 
first discovered by Kantorovich and valid for a large family of cost functions.\ 
The dual variational problem 
permits, among other things, 
to describe the geometry of the optimizers of the original variational
problem \eqref{eq:supell}.

A set $\Gamma \subset X^2$ is cyclically monotone with respect to a cost function $c$ if 
it is $c$-optimal with respect to perturbation by finitely many points (for the precise definition we refer to \cite[Sec. 5]{Vil}).\
Kantorovich duality gives a recipe to construct such cyclically monotone sets.\
The relevance of this condition can be understood by mentioning that, 
in the smooth Riemannian setting, cyclical monotonicity of a transport plan is equivalent to its
optimality.

However, when the cost function is not real valued, 
for instance in our setting with $\ell(x,y)$ as in \eqref{eq:defell}, 
more attention is required.\
In particular, Kantorovich duality is more subtle (\cite{BertrandPratelliPuel}, \cite{KellSuhr}, \cite{McCann}) and the equivalence between 
 the optimality of a transport plan and the cyclical monotonicity of its support breaks down.

To have a full duality theory we have to consider a more regular family of measures.\
The following definition relaxes the notion of $q$-separated introduced by McCann 
\cite[Definition 4.1]{McCann} in the smooth Lorentzian setting.

\begin{definition}[Timelike $p$-dualisable]\label{D:dualisable}
 Let  $(X,\sfd, \ll, \leq, \tau)$ be a Lorentzian pre-length space and 
 let $p\in (0,1]$.\ We say that $(\mu,\nu)\in \mathcal{P}(X)^{2}$ 
 is \emph{timelike $p$-dualisable (by $\pi\in \Pi_{\ll}(\mu,\nu)$)}  
 if 
\begin{enumerate}
\item  $\ell_{p}(\mu,\nu)\in (0,\infty)$;
\item  $\pi\in  \Pi_{\ll}^{p\text{-opt}}(\mu,\nu)$;
\item there exist measurable functions $a,b:X\to \R$, 
with $a\oplus b \in L^{1}(\mu\otimes \nu)$ such that 
$\ell^{p}\leq a\oplus b$ on $\supp \, \mu \times  \supp \, \nu $.
\end{enumerate}
\end{definition}

If $X$ is globally hyperbolic, $(\mu,\nu)$ have compact support and admit a timelike $p$-optimal plan $\pi\in \Pi^{p\text{-opt}}_{\ll}(\mu, \nu)$,  then $(\mu,\nu)$ is timelike 
$p$-dualisable by $\pi$.\ In particular, the notion of timelike $p$-dualisability relaxes the condition of $q$-separation introduced by McCann \cite{McCann}.\
Moreover timelike $p$-dualisabily ensures a weak form of duality (for the proof see \cite[Prop. 2.19]{CaMo:20}):

\begin{proposition}[Weak Kantorovich duality] 
Fix $p\in (0,1]$. Let  $(X,\sfd, \ll, \leq, \tau)$ be a (resp. locally)  causally closed globally hyperbolic Lorentz geodesic space.\ If $(\mu,\nu)\in \mathcal{P}(X)^{2}$ (resp.  $\mathcal{P}_{c}(X)^{2}$) is timelike $p$-dualisable, then Kantorovich duality holds:
 \begin{equation}\label{eq:KantDual}
 \ell_{p}(\mu,\nu)^{p}=\inf \left\{ \int_{X} u\, \mu+ \int_{X} v\, \nu \right\} ,
 \end{equation}
 where the $\inf$ is taken over all measurable functions $u:\supp\, \mu\to \R\cup \{+\infty\}$ and  $v:\supp\, \nu\to \R\cup \{+\infty\}$ with $u\oplus v\geq \ell^{p}$ on $\supp \, \mu \times  \supp \, \nu $ and  $u\oplus v\in L^{1}(\mu\otimes \nu)$.\ Furthermore, the value of the right hand side does not change if one restricts the $\inf$ to bounded and continuous functions. 
 \end{proposition}

We next discuss the validity of the  \emph{strong} Kantorovich duality, i.e. the existence of optimal functions (called Kantorovich potentials)
achieving 
the infimum on the right hand side of   \eqref{eq:KantDual}.\ 

Since the dual minimization problem runs over the couples of functions 
verifying $u\oplus v\geq \ell^{p}$ 
on $\supp \, \mu \times  \supp \, \nu$, one can consider 
special couples of functions.\ The following definition is indeed taylored
to this duality principle.\

\begin{definition}[$\ell^{p}$-concave functions, $\ell^{p}$-transform and  $\ell^{p}$-subdifferential]
Fix $p\in (0,1]$ and let $U,V\subset X$.\  A measurable function $\varphi:U\to \R$  is \emph{$\ell^{p}$-concave relatively to $(U,V)$} if there exists a function $\psi:V \to \R$ such that
$$
\varphi(x)=\inf_{y\in V} \psi(y)- \ell^{p}(x,y), \quad \text{for all } x\in U.
$$
The function 
\begin{equation}\label{eq:defellptransform}
\varphi^{(\ell^{p})}: V\to \R \cup \{-\infty\}, \quad \varphi^{(\ell^{p})}(y):=\sup_{x\in U} \varphi(x)+ \ell^{p} (x,y)
\end{equation}
is called \emph{$\ell^{p}$-transform of $\varphi$}.\ The \emph{$\ell^{p}$-subdifferential $\partial_{\ell^{p}} \varphi\subset (\supp\, \mu\times \supp\, \nu )\cap X^{2}_{\leq}$} is defined by
$$
\partial_{\ell^{p}} \varphi := \{(x,y)\in (\supp\, \mu\times \supp\, \nu )\cap X^{2}_{\leq}\,:\,  \varphi^{(\ell^{p})}(y) - \varphi(x) = \ell^{p}(x,y)\}.
$$
\end{definition}

Then the strong form of duality can be defined as follows.

\begin{definition}[Strong Kantorovich duality]\label{def:StrongKantDual}
Fix $p\in (0,1]$. We say that $(\mu,\nu)\in \mathcal{P}(X)^{2}$ satisfies \emph{strong $\ell^{p}$-Kantorovich duality}  if 
\begin{enumerate}
\item  $\ell_{p}(\mu,\nu)\in (0,\infty)$;
\item  there exists $\varphi: \supp \, \mu\to \R$  which is  $\ell^{p}$-concave  relatively to $(\supp \, \mu, \supp \, \nu)$ and
 satisfying
$$
 \ell_{p}(\mu,\nu)^{p}= \int_{X} \varphi^{(\ell^{p})}(y) \, \nu ({\rm d}y) - \int_{X} \varphi(x) \, \mu({\rm d}x).
$$
\end{enumerate}
\end{definition}

It is immediate to check that if  $(\mu,\nu)\in \mathcal{P}(X)^{2}$  satisfies strong $\ell^{p}$-Kantorovich duality, then the following holds: a plan $\pi\in \Pi_{\leq}(\mu,\nu)$ is $\ell_{p}$-optimal if and only if
\begin{equation*}
\varphi^{(\ell^{p})}(y)  - \varphi(x)  =\ell^{p}(x,y)=\tau(x,y)^{p}, \quad \text{for $\pi$-a.e. $(x,y)$}, 
\end{equation*}
i.e. if and only if $\pi (\partial_{\ell^{p}} \varphi)=1$. 

Motivated by this remark, we have devised a set of couples 
of probability measures larger
than  those  satisfying the strong Kantorovich duality (but still sufficiently regular).

\begin{definition}[Strongly timelike  $p$-dualisable] \label{def:TStrongDual}
A pair $(\mu,\nu) \in (\mathcal{P}(X))^{2}$ is said to be \emph{strongly timelike $p$-dualisable} if
\begin{enumerate}
\item  $(\mu,\nu)$ is timelike $p$-dualisable;
\item  there exists a measurable $\ell^{p}$-cyclically monotone set $\Gamma\subset X^{2}_{\ll} \cap (\supp \, \mu \times \supp \,\nu)$ such that a plan  $\pi\in \Pi_{\leq}(\mu,\nu)$ is $\ell_{p}$-optimal if  and only if $\pi$ is concentrated on $\Gamma$, i.e. $\pi(\Gamma)=1$.
\end{enumerate}
\end{definition}

The next two propositions 
(for their proof see \cite[Cor. 2.29, Cor.  2.30]{CaMo:20}) 
show that the notion of strongly timelike $p$-dualisable measures is non-empty: the first one addresses the ``local'' behaviour. 

\begin{proposition}\label{lem:q-dualPcX}
Fix $p\in (0,1]$. Let  $(X,\sfd, \ll, \leq, \tau)$ be a  causally closed  (resp. locally causally closed) globally hyperbolic Lorentzian geodesic space and assume that $\mu,\nu\in \mathcal{P}(X)$ (resp.  $\mathcal{P}_{c}(X)$) satisfy:
\begin{enumerate}
\item there exist measurable functions $a,b:X\to \R$  with $a\oplus b \in L^{1}(\mu\otimes \nu)$ such that $\tau^{p}\leq a \oplus b$ on $\supp \, \mu \times  \supp \, \nu $;
\item $\supp \, \mu \times \supp \, \nu \subset X^{2}_{\ll}$.
\end{enumerate}
Then $(\mu,\nu)$ is strongly timelike $p$-dualisable.
\end{proposition}

The second one shows that in the case where $\nu$ is a Dirac measure,   strongly timelike $p$-dualisability is  equivalent to timelike  $p$-dualisability.

\begin{proposition}\label{cor:STDualDelta}
Let  $(X,\sfd, \ll, \leq, \tau)$ be a Lorentzian pre-length space and let $p\in (0,1]$.\ Fix $\bar x\in X$ and let $\nu:=\delta_{\bar x}$.\ Assume that $\mu \in \mathcal{P}(X)$ satisfies:
\begin{equation}\label{eq:DualnuDirac}
\tau(\cdot, \bar x)^{p} \in L^{1}(X,\mu) \quad \text{and} \quad \tau(\cdot, \bar{x})>0 \text{ $\mu$-a.e. }. 
\end{equation}
Then $(\mu,\nu)$ is strongly timelike  $p$-dualisable.\ In other words, in the case where $\nu$ is a Dirac measure,  strongly timelike  $p$-dualisability is  equivalent to  timelike  $p$-dualisability.
\end{proposition}

\subsection{Geodesic structure of the Lorentz-Wasserstein space} 

A relevant object in the study of the geometry of Lorentz-Wasserstein space 
are geodesics.

\begin{definition}
Let $(X,\sfd,\ll,\leq,\tau)$ be a Lorentzian pre-length space.
\\We say that 
$(\mu_{s})_{s\in [0,1]}\subset\mathcal{P}(X)$ is an $\ell_{p}$-geodesic if and only if
$$
\ell_{p}(\mu_{s},\mu_{t}) = (t-s)\ell_{p}(\mu_{0},\mu_{1}), \quad \text{ for all } s,t\in [0,1] \text{ with }  s \leq t.
$$
In particular $\ell_{p}$-geodesics are implicitly causal future-directed.
\end{definition}

In the next proposition (for the proof see \cite[Prop. 2.32]{CaMo:20}) we collect some useful properties of $\ell_{p}$-geodesics.\ Before stating it, we introduce the evaluation map 
\begin{equation}\label{def:eet}
\ee_{t}: C([0,1], X) \to X, \quad \gamma\mapsto \ee_{t}(\gamma):=\gamma_{t}, \quad \text{ for all } t\in [0,1],
\end{equation}
and  the stretching/restriction operator ${\rm restr}_{s_{1}}^{s_{2}}: C([0,1], X) \to C([0,1], X)$
\begin{equation}\label{def:restr}
({\rm restr}_{s_{1}}^{s_{2}} \gamma)_{t}:= \gamma_{(1-t)s_{1}+t s_{2}}, \quad \text{ for all } s_{1},s_{2}\in [0,1], s_{1}<s_{2} \text{ and all }t\in [0,1].
\end{equation}

\begin{proposition}\label{prop:ellqgeodCS}
Let  $(X,\sfd, \ll, \leq, \tau)$ be a  $\cK$-globally hyperbolic, Lorentzian geodesic space.\ Let $\mu_{0},\mu_{1}\in \Prob_{c}(X)$ such that there exists $\pi\in \Pi^{p\text{-opt}}_{\leq} (\mu_{0}, \mu_{1})$ with $\supp \, \pi \Subset \{\tau>0\}$ (in particular, if  $\supp (\mu_{0}\otimes \mu_{1})\Subset \{\tau>0\}$).\ Then
\begin{enumerate}
\item There always exists an $\ell_{p}$-geodesic from $\mu_{0}$ to $\mu_{1}$.

\item For every $\ell_{p}$-geodesic $(\mu_{t})_{t\in [0,1]}$ from $\mu_{0}$ to $\mu_{1}$ there exists a probability measure $\eta\in \Prob(C([0,1], X)$ such that  $(\ee_{t})_{\sharp} \eta= \mu_{t}$ for every $t\in [0,1]$ and $\eta$-a.e. $\gamma$ is a maximal causal curve from $\gamma_{0}\in \supp \, \mu_{0}$ to  $\gamma_{1}\in \supp \, \mu_{1}$.\ 
Such an $\eta$ is called \emph{$\ell_{p}$-dynamical optimal plan} and the set of dynamical optimal plans from $\mu_{0}$ to $\mu_{1}$ is denoted by ${\rm OptGeo}_{\ell_{p}}(\mu_{0}, \mu_{1}) $.

\item If $\eta\in {\rm OptGeo}_{\ell_{p}}(\mu_{0}, \mu_{1})$ then  for all  $s_{1},s_{2}\in [0,1]$ with $s_{1}<s_{2}$ it holds:
$$
\eta^{s_{1},s_{2}}:= ({\rm restr}_{s_{1}}^{s_{2}})_{\sharp} \eta\in {\rm OptGeo}_{\ell_{p}}((\ee_{s_{1}})_{\sharp}\eta, (\ee_{s_{2}})_{\sharp}\eta).
$$ 
\end{enumerate}
\end{proposition}
Proposition \ref{prop:ellqgeodCS} proves that the
$\ell_p$-optimal transport is performed along geodesics of the underlying space $X$.\
This link with the geometry of $X$ has been crucially used to study Ricci curvature.

\subsection{Timelike Ricci curvature lower bounds: smooth setting}\label{sec:TRLBsmooth}

In the sequel, we fix a non-negative Borel measure $\mm$ on $(X,\sfd)$ which is finite on bounded sets.\ Such an $\mm$ will play the role of reference volume measure.\ In case of a smooth Lorentzian manifold $(M,g)$, a natural choice is given by $\mm={\rm vol}_g$ (i.e. the volume measure associated to $g$) or $\mm=\exp(f)\; {\rm vol}_g$, where $f\in C^\infty(M)$ plays the role of a weight.\ Once a reference volume measure $\mm$ is fixed, one can define an entropy functional.

\begin{definition}
Given a probability measure $\mu \in \mathcal{P}(X)$ we define 
its relative Boltzmann-Shannon entropy by
\begin{equation}\label{eq:defEnt}
\Ent(\mu|\mm) = \int_{M} \rho \log(\rho) \, \mm,
\end{equation}
if $\mu = \rho \, \mm$ is absolutely continuous with respect to $\mm$ and $(\rho\log(\rho))_{+}$ is  $\mm$-integrable.\ 
Otherwise we set $\Ent(\mu|\mm) = +\infty$.
\end{definition}

A simple application of Jensen's inequality using the convexity of $(0,\infty)\ni t\mapsto t \log t$ gives 
\begin{equation}\label{E:jensenE}
\Ent(\mu|\mm)\geq -\log \mm(\supp \, \mu)>-\infty,\quad \text{for all } \mu\in  \mathcal{P}_{c}(X).
\end{equation}
We set  $\Dom(\Ent(\cdot|\mm)):=\{\mu\in \mathcal{P}(X)\,:\, \Ent(\mu|\mm)\in \R\}$ to be the finiteness domain of the entropy.

Recall that a sequence of probability measures $(\mu_n)\subset \mathcal{P}(X)$ is \emph{narrowly} convergent to $\mu\in \mathcal{P}(X)$ as $n \to \infty$ if
$$
\lim_{n\to \infty} \int_X f \, \mu_n=  \int_X f \, \mu, \quad \text{for all $f\in C^0_b(X)$,}
$$
where $C^0_b(X)$ denotes the space of continuous and bounded real functions defined on $X$.

An important property of the relative entropy is the lower-semicontinuity under narrow convergence,  (for the proof see for instance \cite[Lemma 4.1]{sturm:I}) 
\begin{equation}\label{eq:Entlsc} 
\mu_{n}\to \mui \; \text{narrowly and } \mm\Big( \bigcup_{n\in \N} \supp\, \mu_{n} \Big) <\infty \quad \Longrightarrow \quad  \liminf_{n\to \infty} \Ent(\mu_{n}|\mm) \geq \Ent(\mui|\mm).
\end{equation}

The connection between Ricci curvature lower bounds and convexity properties of the entropy functional along Wasserstein geodesics (classically called “displacement convexity") is now well understood in the Riemannian setting.\ It was proved by McCann (in his Ph.D. thesis \cite{McCannThesis}) that displacement convexity holds in $\R^n$.\ Then Otto-Villani \cite{OV} formally observed that non-negative Ricci implies displacement convexity  for smooth Riemannian manifolds, an observation then rigorously proved by Cordero Erausquin-McCann-Schmuckenschl\"ager \cite{CEMCS}.\ The circle was then closed by von Renesse-Sturm \cite{VrSt}, who proved a complete characterization of Ricci lower bounds for smooth Riemannian manifolds in terms of convexity properties of the entropy.

In the same spirit as in the Riemannian setting,
the convexity properties of the entropy functional along $\ell_p$-geodesics  are equivalent to lower bounds on the Ricci curvature. 
This was proved first by McCann \cite{McCann} for non-negative timelike Ricci lower bounds; in an independent slightly subsequent work, M.-Suhr \cite{MoSu} characterised general lower and upper timelike Ricci bounds (for the proof of the exact statement below, see \cite[Theorem 3.1]{CaMo:20}).

\begin{theorem}[McCann and M.-Suhr]
\label{T:smoothTCD}
\vspace{-0.2cm}
Let $(M^{n},g)$ be a smooth globally hyperbolic spacetime, $0< p <1$ and denote by 
${\rm vol}_g$ the associated volume measure. 
\\Then the following assertions are equivalent:
\begin{itemize}
\item  $\Ric_{g}(v,v) \geq -K g(v,v) $, for every timelike $v \in TM$.

\item $\forall \,(\mu_{0},\mu_{1})\in ({\rm Dom}(\Ent(\cdot|{\rm vol}_{g})))^{2}$ 
timelike $p$-dualisable, there exists a (unique) $\ell_{p}$-geodesic $(\mu_{t})_{t\in [0,1]}$  s.t. $[0,1]\ni t\mapsto e(t) : = \Ent(\mu_{t}|{\rm vol}_{g})$ is semi-convex 
(hence locally Lipschitz on $(0,1)$) and it satisfies:
\begin{equation}\label{eq:conve}
e''(t) - \frac{1}{n} e'(t)^{2 } \geq K \int_{M\times M} \tau(x,y)^{2} \, \pi({\rm d}x{\rm d}y),
\end{equation}
in the distributional sense on $[0,1]$.
\item For any couple $(\mu_{0},\mu_{1})\in ({\rm Dom}(\Ent(\cdot|{\rm vol}_{g}))\cap \Prob_{c}(X))^{2}$ which is strongly timelike $p$-dualisable there exists a (unique) $\ell_{p}$-geodesic $(\mu_{t})_{t\in [0,1]}$ joining them and satisfying \eqref{eq:conve}.
\end{itemize} 
\end{theorem}

\begin{remark}[Timelike Ricci lower bounds and energy conditions in General Relativity]
$\phantom{as}$

\begin{itemize}
\item \textbf{Weak energy condition}. A lower bound on the timelike Ricci curvature of a spacetime $(M^{n},g)$, i.e. the first item in Theorem \ref{T:smoothTCD}:
\begin{equation}\label{eq:Ricgeq-Kg}
\text{There exists $K\in \R$ such that } \Ric_{g}\geq- K g(v,v)  \text{ for all  timelike vectors $v\in TM$,}
\end{equation}
is quite a natural assumption in General Relativity.\ Of course, for a $C^{2}$-metric $g$,  the lower bound \eqref{eq:Ricgeq-Kg} is satisfied on compact subsets of the space-time.
 
Recalling that the Einstein equations postulate proportionality of $\Ric_{g}$ and $T-\frac{1}{n-2} \tr_{g}(T) g$ (where $T$ is the so-called energy-momentum tensor),  for a general cosmological constant $\Lambda \in \R$, the lower bound \eqref{eq:Ricgeq-Kg} is equivalent to requiring that 
 $$T(v,v) \geq  - \frac{1}{n-2} \tr_{g}(T) + \frac{1}{8\pi} \left(  K-\frac{2\Lambda}{n-2}\right),  \text{ for all  $v\in TM$ with $g(v,v)=-1$}.$$
 In particular, if  $\inf_{M} \tr_{g}(T)>-\infty$ (or, equivalently,  $\inf_{M} {\rm R}_{g}>-\infty$ where  ${\rm R}_{g}$ is the scalar curvature of $g$), then the \emph{weak} energy condition $T(v,v)\geq 0$ for all timelike $v$ (which is believed to hold for most physically reasonable $T$, according to \cite[pag. 218]{Wald}) implies \eqref{eq:Ricgeq-Kg}.
 
\item \textbf{Hawking-Penrose's strong energy condition} (SEC for short).\ The SEC  asserts that, calling $T$ the energy-momentum tensor in the Einstein equations, it holds $T(v,v)\geq \frac{1}{2}\tr_{g}(T)$ for every time-like vector $v\in TM$  satisfying $g(v,v)=-1$.\  Assuming that the space-time $(M,g)$ satisfies the Einstein equations  with zero cosmological constant, the SEC is equivalent to $\Ric_{g}(v,v)\geq 0$ for every timelike vector $v\in TM$.\ This corresponds to the case $K=0$ in Theorem \ref{T:smoothTCD}.

The SEC, proposed  by Hawking and Penrose  \cite{Pen,Haw66, HawPen70},  plays a key role in General Relativity.\ For instance, in the presence of trapped surfaces, it  implies that the space-time is  geodesically incomplete.\ This fact is interpreted as singular behaviour possibly connected to the presence of a black hole (for a general discussion about singularity theorems see for instance the monographs \cite{HawEll,Wald, La}).\ Moreover, due to the averaged-focusing effect on geodesics, the SEC is sometimes interpreted as the geometric counterpart of the fact that gravity is an attractive force.
\end{itemize}

\end{remark}

\begin{remark}[A heuristic thermodynamic interpretation of Theorem \ref{T:smoothTCD}]
%
An $\ell_{p}$-geodesic  
$(\mu_{t})_{t\in [0,1]}$  can be interpreted as  
the evolution \footnote{strictly speaking $t$ is not the proper time, but only a variable parametrising the evolution} 
of a distribution (in space-time) of a gas.\
Theorem \ref{T:smoothTCD} says that timelike Ricci lower bounds 
(which correspond to  energy conditions in General Relativity) 
can be equivalently formulated in terms of the convexity properties of the 
Boltzmann-Shannon entropy along such 
evolutions $(\mu_{t})_{t\in [0,1]}\subset \Prob(M)$.\ 
Extrapolating a bit more,  we might say that the second law of 
thermodynamics (i.e. in a natural thermodynamic process, the sum of the entropies of the interacting thermodynamic systems decreases, due to our sign convention) 
concerns the \emph{first} derivative of the Boltzmann-Shannon entropy; 
gravitation (in the form of Ricci curvature) is instead related to the \emph{second} order derivative of the Boltzmann-Shannon entropy along a natural thermodynamic process.
\end{remark}

\begin{remark}[Some related physics literature]

The existence of strong connections between thermodynamics and General Relativity is not new in the physics literature; it has its origins at least in the work Bekenstein \cite{Bek} and Hawking with collaborators \cite{Haw} in the mid-1970s  about black hole thermodynamics.\  These works inspired a new research field in theoretical physics, called entropic gravity (also known as emergent gravity), asserting  that gravity is an entropic force rather than a fundamental interaction.\ 
Let us give a  brief account.\  In 1995  Jacobson \cite{Jac} derived the Einstein equations from the proportionality of entropy and horizon area of a black hole, exploiting the fundamental relation $\delta Q=T \, \delta S$  linking heat $Q$, temperature $T$ and entropy $S$.\  Subsequently, other physicists (let us mention Padmanabhan \cite{Pad}) have been exploring links between gravity and entropy.\
More recently, in 2011 Verlinde \cite{Ver} proposed a heuristic argument suggesting that (Newtonian) gravity can be identified with an entropic force caused by changes in the information associated with the positions of material bodies.\ A
relativistic generalization of those arguments  leads to the Einstein equations.

Theorem \ref{T:smoothTCD} can be seen as an additional strong connection between  general relativity and thermodynamics/information theory.\ It would be interesting to explore this relationship further.
\end{remark}

\begin{remark}
Let us stress that the framework \eqref{eq:Ricgeq-Kg} includes any solution of the Einstein equations \emph{in vacuum} (i.e. with null stress-energy tensor $T$) with \emph{possibly non-zero cosmological constant} $\Lambda$.\ Already such a framework is highly interesting as the standard black hole metrics (e.g. Schwartzshild, Kerr) are solutions of the Einstein vacuum equations, and also the more recent literature on black holes typically focuses on vacuum solutions (see e.g. \cite{Chri, Daf, DHR, KRS}).\  A key role in such  breakthroughs on black holes  is given by a deep analysis of the system of non-linear hyperbolic partial differential equations corresponding to the Einstein vacuum equations (in a suitable gauge).\ At least in the smooth setting, it was recently proved by the second author and Suhr \cite{MoSu} that the optimal transport point of view is compatible with the hyperbolic PDEs one (in the sense that it is possible to characterise solutions of the Einstein equations in terms of optimal transport).\  We believe that all of this suggests that the proposed optimal transport approach has potential for the future.
\end{remark}

\section{Synthetic Theory of timelike Ricci lower bounds}\label{Sec:Ricci>Synt}

\subsection{Definitions of $\mathsf{TCD}^{e}_{p}(K,N)$,  $\mathsf{wTCD}^{e}_{p}(K,N)$ and $\TMCP^{e}(K,N)$}

Theorem \ref{T:smoothTCD} proves the equivalence between a condition requiring the 
smoothness of the Lorentzian metric $g$ and a 
condition (the ordinary differential inequality for $e$) 
that can be formulated in any Lorentzian pre-length space.\ 
Inspired by the Riemannian counterpart which led to the flourishing theory of Lott-Sturm-Villani $\CD(K,N)$ metric measure spaces \cite{sturm:I, sturm:II, lottvillani:metric},  it is then natural to turn it into a definition.

\begin{definition}[$\TCD$ condition]\label{def:TCD(KN)}
Fix $p\in (0,1)$, $K\in \R$, $N\in (0,\infty)$.\  
A  measured pre-length space 
$(X,\sfd,\mm, \ll, \leq, \tau)$ satisfies the $\mathsf{TCD}^{e}_{p}(K,N)$ 
(resp. $\mathsf{wTCD}^{e}_{p}(K,N)$)  if 
for any couple $(\mu_{0},\mu_{1})\in {\rm Dom}(\Ent(\cdot|\mm))^{2}$ which is 
timelike $p$-dualisable   (resp. $(\mu_{0},\mu_{1})\in 
({\rm Dom}(\Ent(\cdot|\mm))\cap  \Prob_{c}(X))^{2}$  
which is   strongly timelike $p$-dualisable) by some 
$\pi\in \Pi^{p\text{-opt}}_{\ll}(\mu_{0},\mu_{1})$,  
there exists an  $\ell_{p}$-geodesic $(\mu_{t})_{t\in [0,1]}$ such that
 $[0,1]\ni t\mapsto e(t) : = \Ent(\mu_{t}|\mm)$ is 
semi-convex (hence locally Lipschitz on $(0,1)$) and it satisfies
\begin{equation}\label{E:differential}
e''(t) - \frac{1}{N} e'(t)^{2 } \geq K \int_{X\times X} \tau(x,y)^{2} \, \pi({\rm d}x{\rm d}y),   \end{equation}
in the distributional sense on $[0,1]$. 
\end{definition}

The condition \eqref{E:differential} is a differential, and therefore infinitesimal,
condition.\ 
One can however formulate an  equivalent global condition taking inspiration from
the 
entropic formulation of the curvature-dimension condition $\CD$ obtained by Erbar-Kuwada-Sturm \cite{EKS}.
This is done by noticing that 
$t\mapsto e(t)$ is semi-convex and  satisfies the inequality 
\eqref{E:differential} if and only if
$t\mapsto u_{N}(t):= \exp(- e(t)/N)$ is semi-convex and satisfies
\begin{equation}\label{eq:sN''}
u_{N}''\leq -\frac{K}{N} \|\tau\|^{2}_{L^{2}(\pi)} \,  u_{N}.
\end{equation} 
Set
\begin{equation}\label{eq:deffsfc}
\fs_{\kappa}(\vartheta):=
\begin{cases}
\frac{1}{\sqrt{\kappa}} \sin(\sqrt{\kappa} \vartheta),  \quad & \kappa>0\\
\vartheta, &\kappa=0\\
\frac{1}{\sqrt{-\kappa}} \sinh(\sqrt{-\kappa} \vartheta), \quad  &\kappa<0\\
\end{cases}, \qquad 
\fc_{\kappa}(\vartheta):=
\begin{cases}
\cos(\sqrt{\kappa} \vartheta),  \quad & \kappa\geq 0\\
\cosh(\sqrt{-\kappa} \vartheta) \quad  &\kappa<0\\
\end{cases}, 
\end{equation}
and
\begin{equation*}
\sigma_{\kappa}^{(t)}(\vartheta):=
\begin{cases}
\frac{\fs_{\kappa}(t\vartheta)}{\fs_{\kappa}(\vartheta)}, \quad & \kappa\vartheta^{2}\neq 0 \text{ and } \kappa\vartheta^{2}<\pi^{2} \\
t,\quad & \kappa \vartheta^{2}=0\\
+\infty \quad &  \kappa \vartheta^{2}\geq \pi^{2}
\end{cases}.
\end{equation*}
Note that the function $\kappa\mapsto \sigma_{\kappa}^{(t)}(\vartheta)$ is non-decreasing for every fixed $\vartheta, t$.\
With the above notation, the differential inequality \eqref{eq:sN''}  is equivalent to the integrated version (cf. \cite[Lemma 2.2]{EKS}):
\begin{equation}\label{eq:sNconc}
u_{N}(t) \geq \sigma^{(1-t)}_{K/N} \left(\|\tau\|_{L^{2}(\pi)}\right) u_{N}(0) + \sigma^{(t)}_{K/N} \left( \|\tau\|_{L^{2}(\pi)} \right) u_{N}(1).
\end{equation} 
We thus proved the following proposition.

\begin{proposition}\label{prop:equivTCD}
Fix $p\in (0,1)$, $K\in \R$ and $N\in (0,\infty)$.\ The measured Lorentzian pre-length space 
$(X,\sfd, \mm, \ll, \leq, \tau)$ satisfies $\mathsf{TCD}^{e}_{p}(K,N)$ 
(resp. $\mathsf{wTCD}^{e}_{p}(K,N)$)
if and only if  for any couple $(\mu_{0},\mu_{1})\in \big(\Dom(\Ent(\cdot|\mm))\big)^{2}$ which is timelike $p$-dualisable  
(resp. $(\mu_{0},\mu_{1})\in [\Dom(\Ent(\cdot|\mm))\cap  \Prob_{c}(X)]^{2}$  which is   
strongly timelike $p$-dualisable) by some   $\pi\in \Pi^{p\text{-opt}}_{\ll}(\mu_{0},\mu_{1})$, 
there exists an $\ell_{p}$-geodesic $(\mu_{t})_{t\in [0,1]}$ such that  the function $[0,1]\ni t\mapsto u_{N}(t) : = U_{N}(\mu_{t}|\mm)$  satisfies \eqref{eq:sNconc}.
\end{proposition}

By considering $(K,N)$-convexity properties only of those $\ell_{p}$-geodesics 
$(\mu_{t})_{t\in [0,1]}$ where $\mu_{1}$ is a Dirac measure one obtains the following weaker condition \cite{CaMo:20} (inspired by the Riemannian counterparts by Sturm \cite{sturm:II} and Ohta \cite{Ohta1}).

\begin{definition}\label{eq:defTMCP}
Fix $p\in (0,1)$, $K\in \R$, $N\in (0,\infty)$.\ The measured globally hyperbolic Lorentz geodesic space $(X,\sfd, \mm, \ll, \leq, \tau)$ satisfies $\mathsf{TMCP}^{e}(K,N)$ if and only if  for any $\mu_{0}\in \Prob_{c}(X)\cap \Dom(\Ent(\cdot|\mm))$ and for any $x_{1}\in X$ such that  $x\ll x_{1}$ for $\mu_{0}$-a.e. $x\in X$, there exists an $\ell_{p}$-geodesic $(\mu_{t})_{t\in [0,1]}$ from $\mu_{0}$ to  $\mu_{1}=\delta_{x_{1}}$
such that 
\begin{equation}\label{eq:defTMCP(KN)}
U_{N}(\mu_{t}|\mm) \geq \sigma^{(1-t)}_{K/N} \left( \sqrt{ \int_{X} \tau(x,x_{1})^{2} \, \mu_{0}({\rm d}x) }\right)\, U_{N}(\mu_{0}|\mm), \quad \text{for all } t\in [0,1).
\end{equation}
\end{definition}

\begin{remark}
The validity of the $\mathsf{TMCP}^{e}(K,N)$ condition is independent of the choice of $p\in (0,1)$ in Definition \ref{eq:defTMCP}. Indeed for any $p,q\in (0,1)$, a curve $(\mu_{t})_{t\in [0,1]}$ with $\mu_{1}=\delta_{\bar{x}}$ is an $\ell_{p}$-geodesic if and only if it is an $\ell_{q}$-geodesic. Let us briefly recall the argument. By Proposition \ref{prop:ellqgeodCS}, given an $\ell_{p}$-geodesic $(\mu_{t})_{t\in [0,1]}$ there exists a measure $\eta\in \Prob( \TGeo(X))$ such that $\mu_{t}=(\mathrm{e}_{t})_{\sharp}\eta$. Thus, for any $q\in (0,1)$, it holds
\begin{align*}
\ell_{q}(\mu_{t}, \mu_{1})^{q}&= \int_{X\times X} \tau^{q}(\gamma_{t}, \bar{x}) \, \eta( {\rm d}\gamma)= (1-t)^{q} \int_{X\times X} \tau^{q}(\gamma_{0}, \bar{x}) \, \eta( {\rm d}\gamma)  \\
&= (1-t)^{q} \int_{X\times X} \tau^{q}(x, \bar{x}) \,\mu_{0}( {\rm d}x)=  (1-t)^{q} \ell_{q}(\mu_{0}, \mu_{1})^{q},
\end{align*}
where in the second identity of the first line we used that $\eta$-a.e. $\gamma$ is a timelike geodesic from $\gamma_{0}$ to $\bar{x}$, and  in the last identity we used that $\mu_{0}\otimes \delta_{\bar{x}}$ is the unique (and thus optimal) plan from $\mu_{0}$ to $\delta_{\bar{x}}$.

\end{remark}

Also the $\TMCP$ condition is able to characterize the lower bounds 
on the time-like Ricci curvature.\ 
In particular (see \cite[Cor. A.2]{CaMo:20}
if $(M^{n}, g)$ is a  globally hyperbolic smooth spacetime of dimension 
$n\geq 2$ without boundary then: 
\begin{enumerate}
\item If $\Ric_{g}(v,v)\geq -K g(v,v)$ for every timelike vector $v\in TM$, then the  associated Lorentzian geodesic space satisfies $\TMCP^{e}(K',N')$ for every $K'\leq K$ and $N'\geq N$.
\item If the  Lorentzian geodesic space associated to $(M^{n}, g)$ satisfies $\TMCP^{e}(K,N)$, then $n\leq N$.
\end{enumerate}

Hence the notion of $\TMCP$ is compatible with the smooth setting.\
Finally, as expected, the $\mathsf{wTCD}^{e}_{p}(K,N)$ condition  implies the 
$\mathsf{TMCP}^{e}(K,N)$ (see \cite[Prop. 3.11]{CaMo:20} for the proof). 

\begin{proposition}\label{prop:CD->MCP}
Fix $p\in (0,1)$, $K\in \R$, $N\in (0,\infty)$.\ The  $\mathsf{wTCD}^{e}_{p}(K,N)$ condition implies the $\mathsf{TMCP}^{e}(K,N)$ condition for locally causally closed, $\sfd$-compatible, $\cK$-globally hyperbolic Lorentzian geodesic spaces.
\end{proposition}

A number of geometric consequences can be deduced for 
measured pre-length spaces verifying the $\wTCD$ 
condition (some of them valid also under the weaker $\mathsf{TMCP}$).\ 
Referring to \cite{CaMo:20} for the complete list, here we only mention 
the timelike  Brunn-Minkowski  inequality, 
the timelike Bishop-Gromov  inequality, 
and the timelike  Bonnet-Myers.

\subsection{Stability of synthetic Ricci curvature lower bounds}\label{sec:stabLB}

One of the main motivations for the huge impact of the Lott-Sturm-Villani theory of synthetic curvature bounds lies in the 
the stability property of the $\CD$ condition: 
if a sequence of metric measure spaces $(X_n,\sfd_n,\mm_n)$ verifies 
$\CD(K,N)$ and it is converging in the measured-Gromov-Hausdorff sense to 
a limit metric measure space $(X,\sfd,\mm)$, then
$(X,\sfd,\mm)$ satisfies the $\CD(K,N)$ condition  as well. 

We have therefore faced the challenging problem of developing a stability 
property for the $\TCD$ condition.\ While there is a well established theory of convergence for Riemannian manifolds (and more generally for metric measure spaces), the situation for Lorentzian manifolds (and even more for Lorentzian length spaces) is more open to investigation.\ 
Without a well established notion of convergence for Lorentzian length spaces already present in the literature, 
we have devised the following one.\ Before stating it, recall that a topological 
embedding is a map $f:X\to Y$ between two topological spaces $X$ and $Y$ such that 
$f$ is continuous, injective and with continuous inverse between $X$ and $f(X)$. We also say that a space $X$ is \emph{pointed}, if a reference point $\star\in X$ is specified.

\begin{definition}\label{D:convergence}
 A sequence of  pointed measured Lorentzian geodesic spaces
$\{(X_{j},\sfd_{j}, \mm_{j}, \star_{j}, \ll_{j}, \leq_{j}, \tau_{j})\}_{j\in \N}$  
converges to
$(X_{\infty},\sfd_{\infty}, \mm_{\infty}, \star_{\infty}, \ll_{\infty}, \leq_{\infty}, \tau_{\infty})$ if and only if:
\begin{enumerate}
    \item There exists a  locally causally closed, $\mathcal{K}$-globally
hyperbolic  Lorentzian geodesic space 
\\ $(\bar X, \bar \sfd,  \overline \ll, \bar \leq, \bar  \tau)$ such that each $j\in \N\cup\{\infty\}$, 
$(X_{j},\sfd_{j}, \mm_{j}, \ll_{j}, \leq_{j}, \tau_{j})$ is 
 isomorphically embedded in it, i.e. there exist 
 topological embedding  maps $\iota_{j}: X_{j}\to \bar{X}$  such that
\begin{itemize}
\item  $x^{1}_{j} \leq_{j} x^{2}_{j}$ 
if and only if $\iota_{j}(x^{1}_{j}) \bar\leq  \iota_{j}(x^{2}_{j})$, for every $j\in \N\cup\{\infty\}$, 
for every   $x^{1}_{j}, x^{2}_{j}\in X_{j}$; 
\item $ \tau_{j} (x^{1}_{j}, x^{2}_{j}) = \bar{\tau}  (\iota_{j}(x^{1}_{j}), \iota_{j}(x^{2}_{j}))$ for every $x^{1}_{j}, x^{2}_{j}\in X_{j}$,  for every $j\in \N\cup\{\infty\}$;
\end{itemize}

\item The measures $(\iota_{j})_{\sharp} \mm_{j}$  converge to $(\iota_{\infty})_{\sharp} \mm_{\infty}$ weakly in duality with $C_{c}(\bar X)$ in $\bar X$, i.e.
$$
\int \varphi\; (\iota_{j})_{\sharp} \mm_{j} \to \int \varphi \; (\iota_{\infty})_{\sharp} \mm_{\infty} \quad \text{for all } \varphi \in C_{c}(\bar X).
$$
\item Convergence of the reference points: $\iota_j(\star_j)\to \iota_\infty(\star_\infty)$ in $\bar{X}$.
\end{enumerate}
\end{definition}

In Definition \ref{D:convergence}
we use topological embeddings to identify
spaces with their image inside a larger space.\ 
This procedure to compare different spaces is rather standard in the framework of metric measure spaces and it provides an effective way to define and study the pointed measured Gromov-Haudorff convergence (in the latter, the assumption is actually that each embedding is \emph{isometric}, which of course is much stronger than merely a topological embedding; for an overview of equivalent formulations of convergence for metric measure spaces see for instance \cite[Sec. 3]{GMS2013}).\ An important assumption is that the causal structures are preserved under the image of the embeddings.

Even though we haven't specifically listed any topological assumption on the sequence of spaces 
$X_{j}$, they actually inherit them from $\bar{X}$ via the topological embeddings $\iota_{j}$.\ 
The map $\iota_{j}$ preserves both the causal relations and $\tau_{j}$, hence
$(X_{j},\sfd_{j}, \mm_{j}, \ll_{j}, \leq_{j}, \tau_{j})$ are 
locally causally closed and $\K$-globally hyperbolic Lorentzian geodesic (by assumption) spaces.

\begin{theorem}[Weak stability of $\mathsf{TCD}$]\label{Thm:WeakStabTCD}
Let  $\{(X_{j},\sfd_{j}, \mm_{j}, \star_{j}, \ll_{j}, \leq_{j}, \tau_{j})\}_{j\in \N}$ be a sequence of pointed measured Lorentzian geodesic spaces 
converging to 
$(X_{\infty},\sfd_{\infty}, \mm_{\infty}, \star_{\infty}, \ll_{\infty}, \leq_{\infty}, \tau_{\infty})$
in the sense of Definition \ref{D:convergence}.\ 
Assume moreover that there exist $p\in (0,1), K\in \R, N\in (0,\infty)$ such that 
$(X_{j},\sfd_{j}, \mm_{j}, \ll_{j}, \leq_{j}, \tau_{j})$ satisfies 
$\mathsf{TCD}^{e}_{p}(K,N)$,  for each $j\in \N$.  

Then the limit space 
$(X_{\infty},\sfd_{\infty}, \mm_{\infty}, \ll_{\infty}, \leq_{\infty}, \tau_{\infty})$
satisfies   $\mathsf{wTCD}^{e}_{p}(K,N)$.
\end{theorem}

Let us comment on the general strategy to obtain a 
stability result for the $\TCD$ condition; this will highlight 
the reasons why full stability is not yet at our disposal, and motivates the presence of the two distinct definitions (weak and strong) of the timelike curvature-dimension condition (see Definition \ref{def:TCD(KN)}).\ For the detailed proof see \cite[Thm. 3.14]{CaMo:20}.
\medskip

\textbf{Some comments on the proof.}
Start by fixing two probability measures $\mu_{0,\infty},\mu_{1,\infty} \in \mathcal{P}_{c}(\bar X)$ both a.c. with respect to $\mm_{\infty}$.
The goal is to construct an $\ell_{p}$-geodesic $(\mu_{t,\infty})_{t\in [0,1]}$ 
such that the map 
 $[0,1]\ni t\mapsto u_{N}(t):= \exp(- \Ent(\mu_{t,\infty}|\mm_{\infty})/N)$ 
 satisfies the inequality \eqref{eq:sNconc}. 

The classical approach in the metric measure space setting 
(see for instance \cite{sturm:I, sturm:II}) would be to use that 
$\mm_{j}\to \mm_{\infty}$ in $W_{2}^{(\bar X,  \bar \sfd)}$ 
(this can be assumed without any loss of generality) and 
to take an optimal plan $\gamma_{j} \in \mathcal{P}(\bar X \times \bar X)$  between 
$\mm_{j}$ and $\mm_{\infty}$.

Then, denoting by $\mu_{i,\infty} = \rho_{i} \mm_{\infty}$, for $i =0,1$,
one defines two sequences of probability measures:
$$
\mu_{i,j} : = (P_{1})_{\sharp} \left(\rho_{i} \gamma_{j}\right), \quad i = 0 ,1. 
\quad 
\Longrightarrow 
\quad
\mu_{i,j} \ll \mm_{j},
$$
also verifying  $\mu_{i,j} \weak \mu_{i,\infty}$ as $j\to \infty$ for $i = 0,1$.

The general scheme (that in this setting breaks down) would be to consider an $\ell_{p}$-geod $(\mu_{t,j})_{t\in[0,1]}$ 
from $\mu_{0,j}$ to $\mu_{1,j}$ given by the $\mathsf{TCD}(K,N)$ condition of $X_{j}$ 
and then to prove (non-trivial but classical) the upper semi-continuity to the limit 
on $u_{N,j}(t)$ and continuity for $u_{N,j}(0)$ and $u_{N,j}(1)$.

In the present setting, the major difficulty is that $\mu_{0,j}$ and $\mu_{1,j}$ may or may not  be causally related; in particular, one cannot invoke the  
$\mathsf{TCD}(K,N)$ condition between $\mu_{0,j}$ and $\mu_{1,j}$.
This is the main reason why we assume that the spaces in the converging sequence satisfy the $\mathsf{TCD}(K,N)$ (i.e. the slightly stronger version of the condition), and we obtain that the limit space satisfies the $\mathsf{wTCD}(K,N)$ (i.e. the slightly weaker version of the condition).

Alternatively, this issue could be circumvented also by imposing some restrictions on the 
family of geodesics along which the concavity of the entropy has to be checked. 
\hfill$\Box$
\\

If either $\mu_{0}$ or $\mu_{1}$ is a Dirac measure, 
then full stability is granted. This is the content of the next result (for a detailed proof see \cite[Thm 3.12]{CaMo:20})

\begin{theorem}[Stability of $\mathsf{TMCP}$]
Let  $\{(X_{j},\sfd_{j}, \mm_{j}, \star_{j}, \ll_{j}, \leq_{j}, \tau_{j})\}_{j\in \N}$ be a sequence of pointed measured Lorentzian geodesic spaces 
converging to 
$(X_{\infty},\sfd_{\infty}, \mm_{\infty}, \star_{\infty}, \ll_{\infty}, \leq_{\infty}, \tau_{\infty})$
in the sense of Definition \ref{D:convergence}. 
Assume moreover that there exist $p\in (0,1), K\in \R, N\in (0,\infty)$ such that 
$(X_{j},\sfd_{j}, \mm_{j}, \ll_{j}, \leq_{j}, \tau_{j})$ satisfies 
$\mathsf{TMCP}^{e}(K,N)$,  for each $j\in \N$.

Then the limit space 
$(X_{\infty},\sfd_{\infty}, \mm_{\infty}, \ll_{\infty}, \leq_{\infty}, \tau_{\infty})$
satisfies   $\mathsf{TMCP}^{e}(K,N)$ as well.
\end{theorem}


\section{Applications: Hawking's singularity theorem and geo\-me\-tric comparison results in a Lorentzian synthetic setting}\label{S:applications}

\subsection{Some useful preliminary concepts}

In this section we will address a generalization to the synthetic setting, 
namely for those Lorentzian geodesic spaces verifying the $\mathsf{TMCP}^{e}(K,N)$ 
condition,
of the Hawking's singularity Theorem. 

To this aim, we will need an extra regularity assumption on the geodesics, namely 
we will require the space to be timelike non-branching (recall the definition \eqref{eq:defTGeo} of $\TGeo(X)$).

\begin{definition}[Timelike non-branching]
A  Lorentzian pre-length space $(X,\sfd, \ll, \leq, \tau)$ is said to be \emph{forward timelike non-branching}  if and only if for any $\gamma^{1},\gamma^{2} \in \TGeo(X)$, it holds:
$$
\exists \;  \bar t\in (0,1) \text{ such that } \ \forall t \in [0, \bar t\,] \quad  \gamma_{ t}^{1} = \gamma_{t}^{2}   
\quad 
\Longrightarrow 
\quad 
\gamma^{1}_{s} = \gamma^{2}_{s}, \quad \forall s \in [0,1].
$$
It is said to be \emph{backward timelike non-branching} if the reversed causal structure is forward timelike non-branching.\ In the case it is both forward and backward timelike non-branching it is said \emph{timelike non-branching}.
\end{definition}

\begin{remark}
By Cauchy-Lipschitz (or Picard-Lindel\"of) Theorem, it is clear that if $(M,g)$ is a space-time whose Christoffel symbols are locally-Lipschitz (e.g. in case $g\in C^{1,1}$) then the associated synthetic structure is timelike non-branching.\ 
It is expected that for spacetimes with a metric of lower regularity (e.g. $g\in C^{1}$ or $g\in C^{0}$) timelike branching may occur. Similar behaviour could occur in closed cone structures 
when the Lorentz-Finsler norm is not strictly convex (see \cite[Remark 2.8]{Min}).\  For a discussion about geodesics in Lorentzian manifolds of low regularity see for instance \cite{SaSt}. 
\end{remark}


The starting point for the Hawking's singularity theorem is to consider
an \emph{achronal} subset of $X$. 

\begin{definition}[Achronal set]
Let $(X,\sfd, \ll, \leq, \tau)$ be a Lorentzian pre-length space.\ A set $V\subset X$ is  \emph{achronal} if $x\not \ll y$, $\forall x,y\in V$. 
\end{definition}

\begin{remark} $V$  is achronal if and only if  $I^{+}(V)\cap I^{-}(V)= \emptyset$.
\end{remark}

We will need to give a meaning to lower bounds on the mean curvature of $V$
and we therefore need some additional regularity.
\\The following compactness property, originally introduced by Galloway \cite{Ga}  in the smooth setting, will play an important role.

\begin{definition}[Future timelike complete (FTC))]
A subset $V\subset X$ is \emph{future timelike complete} (FTC for short) if for each point  $x\in I^{+}(V)$, the intersection
$J^{-}(x)\cap V \subset V$ has compact closure (w.r.t. $\sfd$) in $V$.\ 
Analogously for \emph{past timelike complete} (PTC).%
\end{definition}

A key role in the proof of the Hawking's singularity theorem is played by the signed time-separation function from an FTC subset $V$.\ This can be seen as the Lorentzian counterpart of the signed distance function from a compact set in a Riemannian manifold.\ Let us recall its definition.

\begin{definition}[Signed time-separation $\tau_V$] 
Let $V\subset X$ be an achronal subset.\ The \emph{signed time-separation} from $V$ is the function  $\tau_{V}:X\to [-\infty, +\infty]$ defined by
\begin{equation}\label{eq:deftauV}
\tau_{V}(x):=
\begin{cases}
\sup_{y\in V} \tau(y,x), &\quad \text{ for }x\in I^{+}(V)\\
-\sup_{y\in V} \tau(x,y),& \quad \text{ for }x\in I^{-}(V) \\
0 &\quad \text{ otherwise}.
\end{cases}
\end{equation}
\end{definition}

Note that $\tau_{V}$ is lower semi-continuous, as supremum of (lower semi-)continuous functions.
\\In order for these suprema to be attained, global hyperbolicity and the geodesic property of $X$ alone are not sufficient; here the additional FTC and PTC assumptions on $V$ are very useful.

\begin{proposition}\label{P:projection}
Let $(X,\sfd, \ll, \leq, \tau)$ be a globally hyperbolic Lorentzian 
geodesic space and 
let $V\subset X$ be an achronal FTC (resp. PTC) subset.\ Then,  for all $x\in I^{+}(V)$ (resp. for all $x\in I^{-}(V)$) there exists 
$y_{x}\in V$ with $\tau_{V}(x)=\tau(y_{x},x)>0$ (resp. $\tau_{V}(x)=-\tau(x, y_{x})<0$). 
\end{proposition}
We refer the reader to \cite[Lemma 1.8]{CaMo:20} for the proof.

\subsection{Disintegration and synthetic mean curvature bounds}

In order to introduce a weak notion of mean curvature, it is useful to define normal variations to the achronal set $V$.\
As a consequence of Proposition \ref{P:projection} and of the reverse triangle inequality \eqref{eq:deftau}, it holds
\begin{equation}\label{eq:tauV1Lip}
\tau_{V}(y)-\tau_{V}(x)\geq \tau(x,y), \quad \text{for all } x,y \in I^{+}(V), \, x\leq y,
\end{equation}
that is  ``$\tau_{V}$ is (reverse) 1-Lipschitz with respect to $\tau$''.\  Relying 
on Proposition \ref{P:projection} one can consider, even in the non-smooth
setting, integral lines $\gamma$ of $\tau_V$ saturating the ``reverse 1-Lipschitz inequality'' \eqref{eq:tauV1Lip}, i.e. timelike geodesics $\gamma$ such that
$$
\tau_{V}(\gamma_{t})-\tau_{V}(\gamma_{s})=\tau(\gamma_{s}, \gamma_{t}), \text{ for all } s\leq t. 
$$
This will provide, up to a set of measure zero, a (non-smooth) “foliation" of $I^+(V)$ by timelike geodesics
that will suffice for our purpose. 

Indeed, removing sets of measure zero (here  the synthetic curvature bound $\mathsf{TMCP}^{e}(K,N)$ is crucially used), 
the following disintegration formula for the reference measure $\mm$ holds true 
(see \cite[Sec. 4]{CaMo:20} for the details):
$$
\mm\llcorner_{I^{+}(V)}
= \int_{Q} h(\alpha,\cdot) \, \mathcal{L}^{1}\llcorner_{X_{\alpha}}\, \qq(d\alpha),
\qquad X_{\alpha} \text{ timelike geodesic for every } \alpha\in Q
$$
\begin{itemize}
\item $Q$ is a (typically uncountable) set of indices, that can be obtained as a Borel quotient set from $X$ under the equivalence relation that identifies points lying on the same integral line of $\tau_V$ saturating the 1-Lipschitz (w.r.t. $\tau$) inequality;
\item $\qq$ is a probability measure over the set of indices $Q$ 
(in case the ambient measure $\mm$ is finite,  this is obtained as push forward of $\mm$ via the aforementioned quotient map after normalization; the case when $\mm$ is not finite can be treated by a cut-off argument);
\item $h(\alpha,\cdot)\in L^{1}_{loc}(X_{\alpha}, \mathcal{L}^{1}\llcorner_{X_{\alpha}})$ for $\qq$-a.e. $\alpha\in Q$;
\item $(X_{\alpha}, |\cdot|, h(\alpha,\cdot)\mathcal{L}^{1})$ verifies $\mathsf{MCP}(K,N)$, i.e. 
\begin{equation}\label{E:MCP0N1d}
\left(\frac{\fs_{K/(N-1)}(b-\tau_{V}(x_{1}))}{\fs_{K/(N-1)}(b-\tau_{V}(x_{0}))}\right)^{N-1}
\leq \frac{h(\alpha, x_{1} ) }{h (\alpha, x_{0})} \leq 
\left( \frac{\fs_{K/(N-1)}(\tau_{V}(x_{1}) - a) }{\fs_{K/(N-1)}(\tau_{V}(x_{0}) -a)} \right)^{N-1},  
\end{equation}
for all $x_{0},x_{1}\in X_{\alpha}$, with $0\leq  a< \tau_{V}(x_{0})<\tau_{V}(x_{1})<b<\pi\sqrt{(N-1)/(K\vee 0)}$.\  
Here $(a,b)\subset \R$ denotes the segment corresponding to the geodesic $X_\alpha$ via the map $\tau_V$, i.e. $(a,b)=\tau_V(X_\alpha)$ (note that the interval can also be closed, or half closed, depending on the geometry of $X$), and $\fs_{K/(N-1)}(\cdot)$ was defined in \eqref{eq:deffsfc}.
\end{itemize}

By a Fubini-Tonelli argument, one can obtain 
the following  coarea-type formula:
\begin{equation}\label{E:disintCoarea}
\mm\llcorner_{I^+_V} = \int_{0}^{\infty} \H_{t} \,{\rm d}t,    
\end{equation}
meaning that for each measurable set $A \subset I^+_V \cup V$ with 
$\mm(A) < \infty$, 
the map $[0, \infty) \ni t \mapsto \H_{t}(A)$ is measurable and 
$$
\mm(A) = \int_{0}^{\infty}\H_{t}(A)\, {\rm d}t  = \int_{0}^{\infty}\H_{t}(A\cap\{\tau_{V} = t\})\, {\rm d}t.
$$ 
Moreover, by construction, 
$\H_{t}$ is concentrated on the level set $\{\tau_{V} = t\}$ 
and  $\H_{0}$ is concentrated on $V$.\
We refer to \cite[Sec. 5.1]{CaMo:20} for the precise definition of the 
family of measures $\H_t$ and all the other missing details. 
\\

We now review how \eqref{E:disintCoarea} is sufficient to 
give a meaning to  
\emph{synthetic mean curvature bounds for $V$} at least in the smooth setting.\ 
In particular we will analyse the disintegration formula and the measures $\cH_{t}$ and their relation to mean curvature bounds.

Let $(M^{n}, g)$ be a $2\leq n$-dimensional smooth globally hyperbolic space-time and  $V\subset M$ be a smooth compact achronal spacelike hypersurface without boundary.\  Then, the signed time-separation function $\tau_{V}$ from $V$  is smooth on a neighbourhood $U$ of $V$ and $\nabla \tau_{V}$ is the smooth timelike past-pointing unit normal vector field along $V$.\ More precisely, 
$$\nabla \tau_{V}(x) \perp T_{x} V, \quad g(\nabla \tau_{V}(x), \nabla \tau_{V}(x))=-1, \quad \text{for all } x\in V.$$
Denote by $\vol_{g}$ the volume measure of $(M^{n},g)$ and by  $\vol_{V}$ the induced $(n-1)$-dimensional volume measure on $V$.\ By compactness of $V$, there exists $\delta>0$ such that the $g$-geodesic $[0, \delta]\ni t\mapsto \exp_{x}(-t \nabla \tau_{V}(x))$ is a future pointing maximal geodesic, for every $x\in V$.\ Define
$${\mathcal U}:=V\times[0, \delta]\subset V\times \R, \quad \Phi:{\mathcal U}\to M, \; \Phi(x,t):=\exp_{x}(-t \nabla \tau_{V}(x)).$$
For $\delta>0$ small enough it is a standard fact (tubular neighbourhood theorem) that $\Phi$ is a diffeomorphism onto its image and that the following integration formula holds true:
\begin{equation}\label{eq:intMvolV}
\int_{M} \varphi\, \vol_{g}=\int_{V}\int_{0}^{\delta} \varphi \circ \Phi(x,t) \, \det D\Phi_{(x,t)}|_{T_{x}V} \, {\rm d}t \, \vol_{V}({\rm d}x), \quad \text{for all } \varphi\in C_{c}(\Phi({\mathcal U})).
\end{equation}
In the smooth setting, we can identify the set of indices (which in general is a Borel quotient set) $Q$ with $V$, and the quotient measure $\qq$ with  $\psi \, \vol_{V}$ defined as follows:
$$
\qq:=  \psi \,\vol_{V} \ll \vol_{V}, \quad \text{where } \psi(x):= \left(\int_{0}^{\delta}   \det D\Phi_{(x,t)}|_{T_{x}V} \, {\rm d}t \right), \,\text{ for all }x\in V.
$$
The integration formula \eqref{eq:intMvolV} can  thus be rewritten as
\begin{equation}\label{eq:intMqq}
\int_{M} \varphi\, \vol_{g}=\int_{V} \frac{1}{\psi(x)} \int_{0}^{\delta} \varphi \circ \Phi(x,t) \, \det D\Phi_{(x,t)}|_{T_{x}V} \, {\rm d}t \, \qq({\rm d}x), \quad \text{for all } \varphi\in C_{c}(\Phi({\mathcal U})).
\end{equation}
By the uniqueness of the disintegration formula, \eqref{eq:intMqq} gives:
$$
h_{\alpha}(t)=\frac{1}{\psi(\alpha)}  \, \det D\Phi_{(\alpha,t)}|_{T_{\alpha}V}, \quad  h_{\alpha}(0)=\frac{1}{\psi(\alpha)}, \quad \text{for all } \alpha\in V \text{and all }t\in [0,\delta].
$$
Moreover, it follows that the measure $\H_{t}$  can be written as
$$
\H_{t} : = \Phi(\cdot, t)_{\sharp} \left( \, \det D\Phi_{(\alpha,t)}|_{T_{\alpha}V}\,  \vol_{V}( {\rm d}\alpha) \right), \quad \text{for all } t\geq 0, 
$$ 
in particular, $\H_{0}=\vol_{V}$ and $\H_{t}$ is the $(n-1)$-volume measure on the hypersurface $\{\Phi(x,t):\, x\in V\}$.

Given any $\phi\in C^{\infty}(V;\R_{\geq 0})$, 
one can consider the region $V_{t,\phi}$ as the domain trapped 
between $V$ and the normal graph of $\phi$.\ 
The first variation of the volume is thus 
$$
\frac{d}{dt} \vol_{g}(V_{t,\phi})=\cH^{n-1}(\{\Phi(x, t\phi(x)):x\in V\}),
$$ 
where  $\cH^{n-1}$ is the standard $(n-1)$-volume of the hypersurface $\{\Phi(x, t\phi(x)):x\in V\}$; in particular,  $\left.\frac{d}{dt}\right|_{t=0} \vol_{g}(V_{t,\phi})=\vol_{V}(V)=\cH_{0}(V)$. 
Taking a further variation, we get 
 \begin{align*}
 \int_{V} \phi^{2}\, g(\vec{H}_{V}, \nabla \tau_{V}) \, \vol_{V} = 
 \left. \frac{d}{dt^{2}}\right|_{t=0}  \vol_{g}(V_{t,\phi})=   
\lim_{t \downarrow 0} \frac{\mm(V_{t,\phi}) - t 
\int_{V} \phi \H_{0}}{t^{2}/2}. 
\end{align*}
 The right hand side, corresponding to the second variation of volume, is thus the first variation of the area which gives the mean curvature $\vec{H}_{V}$ of $V$.\
The right hand side can therefore be chosen as a starting point 
for defining bounds on the mean curvature.\
Still we need to make precise in the non-smooth setting the definition of variation set $V_{t,\phi}$.\ 
We will use the ``initial-point projection map'' 
$\fa: I^+_{V}\to V$ assigning to each $x \in I^+_V$ the point 
given by Proposition \ref{P:projection}.\
It is not hard to check it is $\mm$-measurable (see \cite{CaMo:20} for details).

\begin{definition}[Synthetic mean curvature bounds]\label{D:MeanCurv}
The Borel achronal FTC subset $V\subset X$ has \emph{forward mean curvature bounded below by 
$H_{0} \in \R$} if $\cH_0$ is a non-negative Radon measure  with $(\fa)_{\sharp} \qq\ll \cH_{0}$ and such that for any normal variation 
$$
V_{t,\phi} : = \{ x \in \T_{V} \colon 0 \leq \tau_{V}(x) \leq t \phi(\fa(x)) \},
$$
the following inequality holds true: 
$$
\limsup_{t \downarrow 0} \frac{\mm(V_{t,\phi}) - t 
\int_{V} \phi \H_{0}}{t^{2}/2} \geq H_{0} \int_{V}\phi^{2} \H_{0}, 
$$
for any bounded Borel function $\phi:V\to [0,\infty)$ with compact support.\
Analogously, $V$ has \emph{forward mean curvature bounded above by  
$H_{0} \in \R$} if  $\cH_0$ is a non-negative Radon measure and for any normal variation 
$V_{t,\phi}$ as above 
the following inequality holds true: 
\begin{equation}\label{eq:2ndvarVol}
\liminf_{t \downarrow 0} \frac{\mm(V_{t,\phi}) - t 
\int_{V} \phi \H_{0}}{t^{2}/2} \leq H_{0} \int_{V}\phi^{2} \H_{0}, 
\end{equation}
for any bounded  Borel  function $\phi:V\to [0,\infty)$ with compact support.
\end{definition}

\begin{remark}[Example of a surface with a conical singularity]
The notion of forward mean curvature bound of Definition \ref{D:MeanCurv} should be compared with the recent related definition proposed by Ketterer \cite{Ket:HK}.\ In the notation of \cite{Ket:HK}, in order to have finite bound $H_0$ one needs an interior and exterior ball condition on $V$ (equivalently, in the smooth setting, to a  local $L^{\infty}$ bound on the full second fundamental form), see \cite[Remark 5.9]{Ket:HK}.\ The notion proposed above in Definition \ref{D:MeanCurv} instead works well even if the set $V$ has corners or conical singularities.
\\
For example, one can check  that the set 
$$
V=\{(x,t)\subset \R^{n,1}: t=\alpha|x|\}, \quad \text{for } \alpha\in (0,1), 
$$
in the $(n+1)$-dimensional Minkowski space-time $\R^{n,1}$ is an achronal topological hypersurface, smooth outside of the origin (where it is Lipschitz) and having forward mean curvature bounded above by  
$H_{0}=0$ in the sense of Definition \ref{D:MeanCurv}.\ Notice that for any compact subset, one could choose the upper bound on the mean curvature to be strictly negative, but such an upper bound approaches zero as $|x|\to \infty$.
\end{remark}

\subsection{The results}

Let us define $D_{H_{0}, K,N}>0$ as follows:
\begin{equation}
D_{H_{0},K,N}:=
\begin{cases}
 \frac{\pi}{2} \sqrt{  \frac{N-1}{K}} & \quad  \text{if $K>0$, $N>1$, $H_{0}=0$}\\
 \sqrt{ \frac{N-1}{K}} \cot^{-1} \left(\frac{-H_{0}}{ \sqrt{K(N-1)}} \right) & \quad  \text{if $K>0$, $N>1$, $H_{0}\in \R\setminus\{0\}$}\\
-\frac{N-1}{H_{0}}  &\quad  \text{if $K=0$, $N>1$, $H_{0}<0$}\\
\sqrt{- \frac{N-1}{K}} \coth^{-1} \left(\frac{-H_{0}}{ \sqrt{-K(N-1)}} \right) &\quad  \text{if $K<0$, $N>1$, $H_{0}<-\sqrt{-K(N-1)}$.}
\end{cases}
\end{equation}

Then the following synthetic version of the Hawking's singularity theorem holds true.

\begin{theorem}[Hawking's singularity theorem for $\TMCP^{e}(K,N)$ spaces]\label{thm:HawkSingSint}
Let  $(X,\sfd, \mm, \ll, \leq, \tau)$ be a timelike non-branching,   locally causally closed, $\cK$-globally hyperbolic, Lorentzian geodesic space satisfying $\mathsf{TMCP}^{e}(K,N)$  for some $p\in (0,1),\, K\in \R,\,  N\in [1,\infty)$  and assume that the causally-reversed structure satisfies the same conditions. 
\\ Let $V\subset X$ be a Borel achronal FTC subset having forward mean curvature bounded above by $H_{0}$ in the sense of Definition \ref{D:MeanCurv}.\  If
\begin{enumerate}
\item $K>0$, $N>1$ and $H_{0}\in \R$, or 
\item $K=0$, $N>1$ and $H_{0}<0$, or
\item $K<0$, $N>1$ and $H_{0}<- \sqrt{-K(N-1)}<0$,
\end{enumerate}
\noindent
then for  every $x\in I^{+}(V)$ it holds that $\tau_{V}(x)\leq D_{H_{0},K,N}$. 
In particular, for every timelike geodesic $\gamma\in \TGeo(X)$ with $\gamma_{0}\in V$, the maximal (on the right) domain of definition is contained in $\big[0, D_{H_{0},K,N}\big]$.\ In case $N=1, \, H_{0}<0$, it holds that $I^{+}(V)=\emptyset$.
\end{theorem}

\textbf{Sketch of the proof}\
The proof of Theorem \ref{thm:HawkSingSint} 
follows the following argument: from the assumption on the mean curvature of 
$V$ one can deduce with some effort some estimates on the logarithmic derivative of the densities of the disintegration formula
$\log h '$.\ 
Then, from the disintegration formula and the localization of the curvature bounds, we obtain that $(X_{\alpha}, |\cdot|, h(\alpha,\cdot)\mathcal{L}^{1})$ verifies $\mathsf{MCP}(K,N)$.\  
Putting together this two facts yields 
bounds on the length of $X_{\alpha}$ and giving the claim. 
 \hfill$\Box$
 \\

For completeness, we also report the direct application of 
Theorem \ref{thm:HawkSingSint} to the case of continuous, 
causally plain metrics.\ Recall that causally plain metrics are a subclass (detected by Chr\'usciel-Grant \cite{CG}) of the space of continuous metrics, introduced to avoid pathological causal behaviours (such as null curves emanating from a point covering a set with non-empty interior, a phenomenon called ``bubbling'' and that may happen for H\"older continuous metrics). 

\begin{corollary}[Hawking's singularity theorem for a $C^{0}$-Lorentzian metric]\label{cor:HawSingLCont}
Let $(M,g)$ be a $2\leq n$-dimensional timelike non-branching, globally
hyperbolic, causally plain spacetime with a   $C^{0}$-Lorentzian metric satisfying $\mathsf{TMCP}^{e}(K,N)$ for some $p\in (0,1), K\in \R, N\in (1,\infty)$.\ Assume the causally-reversed structure to 
satisfy the same conditions. 
\\ 
Let $V\subset M$ be a Borel achronal FTC subset (or, more strongly, let $V$ be a Cauchy hypersurface) having forward mean curvature bounded above by $H_{0}<0$ (in the synthetic sense).\ If
\begin{enumerate}
\item $K>0$, $N>1$ and $H_{0}\in \R$, or 
\item $K=0$, $N>1$ and $H_{0}<0$, or
\item $K<0$, $N>1$ and $H_{0}<-\sqrt{-K(N-1)}<0$,
\end{enumerate}
then for  every $x\in I^{+}(V)$ it holds $\tau_{V}(x)\leq D_{H_{0},K,N}$.
In particular, for every timelike geodesic $\gamma\in \Geo(M)$ with $\gamma_{0}\in V$, the maximal (on the right) domain of definition is contained in 
$[0, D_{H_{0},K,N}]$.\ Hence $(M,g)$ is not timelike geodesically complete.
\end{corollary}

\begin{remark}[Literature about Hawking's singularity theorem]\label{rem:litHS}
Hawking's singularity theorem was established in \cite[Thm. 4, p. 272]{HawEll} for smooth space-times (the proof works for $C^{2}$ metrics) assuming that $V$ is a \emph{compact} spacelike slice.\ The result was extended to  $C^{1,1}$ metrics in \cite{Kunz:Haw}  and to $C^1$ metrics in \cite{Graf}, by approximating the metric of low regularity with smoother metrics.\ The extension to non-compact \emph{future causally complete} $V$ was achieved in \cite[Thm. 3.1]{Ga} (see also \cite{GrTr}) in the smooth setting, and extended to $C^{1,1}$ metrics in \cite{Graf1}.\ Theorem  \ref{thm:HawkSingSint} and Corollary \ref{cor:HawSingLCont},  already in the smooth setting, relax the future causal completeness with the weaker future timelike completeness (in addition to extending the results to a  synthetic framework, including $C^{0}$ metrics). 
The Hawking (as well as the Penrose and Hawking-Penrose) singularity theorem was also extended to (smooth) closed cone structures \cite{Min} and smooth weighted Lorentz-Finsler manifolds \cite{LMO}.

A first synthetic singularity theorem was recently shown in \cite{AGKS} under the stronger assumptions that the space is a synthetic warped product with lower bounds on the sectional curvature in the sense of comparison triangles (\'a la Alexandrov).

Few weeks after we announced our work \cite{CaMo:20}, we learnt of \cite{BKMW}, proving a \emph{Riemannian} version of Hawking's singularity theorem in the framework of \emph{metric measure spaces} with Ricci curvature bounded below in synthetic sense via optimal transport.
\end{remark}

\begin{remark}[Recovering the classical Hawking's singularity theorem in the smooth setting] \qquad
\begin{itemize}
\item \emph{Causality condition}. A smooth globally hyperbolic space-time (and its time-reverse structure) enters in the framework of  timelike non-branching, locally causally closed, $\cK$-globally hyperbolic, Lorentzian geodesic spaces (see Section \ref{sec:LSS}).
\\Moreover, if $V$ is a compact achronal spacelike hypersuface in classical terms then it is achronal and future timelike complete in the sense of the present survey.
\item \emph{Strong energy condition}. The classical strong energy conditions amounts to requiring that the timelike Ricci curvature is non-negative, which in turns implies that $\mathsf{TMCP}^{e}(0,N)$ is satisfied both by the space and by its time-reverse structure.
\item \emph{Positive initial expansion of timelike geodesics emanating from $V$}. This condition amounts to requiring that the mean curvature of $V$ in the future pointing direction is positive. Our synthetic notion of positive mean curvature recovers the classical one in the smooth setting (see discussion before Definition \ref{D:MeanCurv}).
\end{itemize}
Combining the discussion above with Theorem \ref{thm:HawkSingSint}, one obtains the classical version of Hawking's singularity theorem for smooth globally hyperbolic spacetimes.
\end{remark}

With similar techniques and ideas one can also obtain some important 
geometric inequalities.\ We report, only in the case of continuous metrics, the Bishop-Gromov inequality and the Bonnet-Myers Theorem, 
the very same statements (as well as other comparison results) 
hold more generally for timelike non-branching, 
locally causally closed, $\cK$-globally hyperbolic, Lorentzian geodesic spaces satisfying $\mathsf{TMCP}^{e}(K,N)$ 
for some $p\in (0,1),\, K\in \R,\,  N\in [1,\infty)$ 
and whose causally-reversed structure satisfies the same conditions (see \cite[Sec. 5.3]{CaMo:20}).

In order to state the results precisely, let us introduce some notation. 
Given $x\leq y\in X$ we set
$$ 
\fI(x,y,t):=\{\gamma_{t}\,:\, \gamma\in \Geo(X),\, \gamma_{0}=x, \, \gamma_{1}=y \}.  
$$
Fix $x_{0}\in X$ and let 
$$
B^{\tau}(x_{0},r):=\{x\in I^{+}(x_{0})\cup \{x_{0}\}: \tau(x_{0}, x)< r\}
$$ 
be the $\tau$-ball of radius $r$ and center $x_{0}$.\   Since typically the volume of a $\tau$-ball is infinite (e.g. in Minkowski space  it is the region below a hyperboloid), it is useful to localise volume estimates using star-shaped sets.\ To this aim, we say that    $E\subset I^{+}(x_{0})\cup\{x_{0}\}$ is $\tau$-star-shaped with respect to $x_{0}$ if 
$ \fI(x_{0},x,t) \subset E$ for every $x\in E$ and $t \in (0,1]$.
Denote by
$$
v(E,r):=\mm(B^{\tau}(x_{0},r)\cap E)
$$
the volume of the $\tau$-ball of radius $r$ intersected with the compact subset  $E\subset I^{+}(x_{0})\cup\{x_{0}\}$, $\tau$-star-shaped with respect to $x_{0}$.

\begin{corollary}[Timelike Bishop-Gromov and Bonner-Myers Theorems]
Let $(M,g)$ be a $2\leq n$-dim.\ timelike non-branching, globally hyperbolic, causally plain spacetime with $C^{0}$-metric satisfying $\mathsf{TMCP}^{e}(K,N)$.\  Assume the causally-reversed structure  
satisfes the same conditions.\ Then the following hold:

\begin{itemize}
\item \emph{Timelike Bishop-Gromov inequality}:  
for every $x_{0}\in M$, every compact subset  $E\subset I^{+}(x_{0})\cup\{x_{0}\}$ $\tau$-star-shaped with respect to $x_{0}$  and every $r<R\leq \pi \sqrt{\frac{N-1}{K \vee 0}}$, it holds
$$
\frac{v(E, r)}{v(E,R)} \geq \frac{\int_{0}^{r} \mathfrak{s}_{K/(N-1)} (t)^{N-1} {\rm d}t}   {\int_{0}^{R}\mathfrak{s}_{K/(N-1)} (t)^{N-1} {\rm d}t }.
$$

\item \emph{Timelike Bonnet-Myers upper bound}:  if $K> 0$, then
$
\sup_{x,y\in X} \tau(x,y) \leq \pi \sqrt{(N-1)/K}.
$
In particular, for any causal curve $\gamma$ it holds 
$\LL_{\tau}(\gamma)\leq \pi \sqrt{(N-1)/K}$ and 
$(M,g)$ is not timelike geodesically complete.
\end{itemize}
\end{corollary}

\section{A synthetic formulation of Einstein's vacuum equations}\label{S:SEE}

The goal of this section is to review the synthetic formulation of Einstein's vacuum equations (i.e. zero stress-energy tensor $T\equiv 0$ but possibly non-zero cosmological constant $\Lambda$) proposed in \cite[App. B]{MoSu}, as a natural outcome of the combination of \cite{MoSu} and \cite{CaMo:20}. 
\\The basic idea is that the Einstein vacuum equations (say with zero cosmological constant for convenience of presentation) can be characterised by the condition 
\begin{equation}\label{eq:Ric0timelike}
\Ric_g(v,v)=0 \qquad \text{for every timelike vector $v\in TM$}. 
\end{equation}
Of course,  \eqref{eq:Ric0timelike} is equivalent to requiring both the timelike Ricci lower bound
\begin{equation}\label{eq:Ric>0timelike}
\Ric_g(v,v)\geq 0 \qquad  \text{for every timelike vector $v\in TM$} 
\end{equation}
and the timelike Ricci upper bound
\begin{equation}\label{eq:Ric<0timelike}
\Ric_g(v,v)\leq 0 \qquad  \text{for every timelike vector $v\in TM$}. 
\end{equation}
The idea is to characterise both bounds separately in a synthetic way, and to ask the validity of both the corresponding synthetic characterizations in the non-smooth setting.

The synthetic  timelike Ricci lower bound \eqref{eq:Ric>0timelike} was already characterised in a synthetic way via optimal transport in Section \ref{Sec:Ricci>Synt}.\ Below we treat the upper bound \eqref{eq:Ric<0timelike}.

\subsection{Synthetic time-like Ricci upper bounds}\label{SS:UB}

The following definition is inspired by Sturm's approach to Riemannian/metric Ricci curvature upper bounds \cite{StUB} and by the  characterization of smooth Lorentzian of timelike Ricci upper bounds obtained in \cite[Thm 4.7, Rem. 4.8]{MoSu}.\  We denote the metric ball in $(X,\sfd)$ with center $x$ and radius $r$ by $B^{\sfd}_{r}(x)$.

 \begin{definition}[Synthetic time-like Ricci upper bounds]\label{def:RUB}
Fix $p\in (0,1)$, $K\in \R$.\ We say that  a  measured  Lorentzian  pre-length space $(X,\sfd,\mm, \ll, \leq, \tau)$ has time-like Ricci curvature bounded above by $K$ in a synthetic sense if there exists $r_{0}>0$ and  a function $\omega:[0,r_{0})\to [0,\infty)$ with $\lim_{r\downarrow 0} \omega(r)=0$ such that for every $r\in [0,r_{0})$ the following holds. 
\begin{itemize}
\item For every $x,y\in X$ with $\sfd(x,y)=r>0$ and such that $B^{\sfd}_{r^4}(x)\times B^{\sfd}_{r^2}(y)\subset X^{2}_{\ll}$,
\item for every $\mu_{0}\in {\rm Dom}(\Ent(\cdot|\mm))$ with $\supp\, \mu_{0}\subset B^{\sfd}_{r^4}(x)$,
\end{itemize}
there exists an $\ell_{p}$-geodesic $(\mu_{t})_{t\in [-1,1]}$ satisfying
\begin{itemize}
\item  $\supp\, \mu_{1}\subset B^{\sfd}_{r^2}(y)$, 
\item $\supp\, \mu_{-1}\times \supp \, \mu_{1}\subset X^{2}_{\leq},$
\item  $\bigcup_{t\in [-1, 1]} \supp \, \mu_{t} \subset B_{10 r_{0}}^{\sfd} (x) $,
\item $
\Ent(\mu_{-1}|\mm)-2 \Ent(\mu_{0}|\mm)+\Ent(\mu_{1}|\mm) \leq (K+\omega(r)) \, r^{2}.
$
\end{itemize}
\end{definition}

A key property of the above notion of time-like Ricci bounded above is the stability under weak convergence of measured  Lorentzian  pre-length spaces.\ The stability of Riemannian/metric Ricci upper bounds via optimal transport was proved by Sturm \cite{StUB}.\ The notion of convergence we use below is a slight reinforcement (we ask that the immersion maps are isometries with respect to the metric structures instead of merely topological embedding maps) of the weak convergence used in Section \ref{sec:stabLB} to show stability of synthetic timelike Ricci \emph{lower} bounds; in any case it is a  natural adaptation to the Lorentzian setting of the mGH convergence used for metric measure spaces (see for instance \cite[Sec. 3]{GMS2013} for an overview of equivalent formulations of convergence for metric measure spaces).

\begin{theorem}[Stability of time-like Ricci curvature upper bounds]\label{thm:stabUB}

Let  $\{(X_{j},\sfd_{j}, \mm_{j}, \star_{j}, \ll_{j}, \leq_{j}, \tau_{j})\}_{j\in \N\cup\{\infty\}}$ be a sequence of pointed measured Lorentzian geodesic spaces  satisfying the following properties :
\begin{enumerate}
\item There exists a  locally causally closed, globally hyperbolic  Lorentzian geodesic space $(X, \sfd,  \ll, \leq,  \tau)$ such that each $(X_{j},\sfd_{j}, \mm_{j}, \ll_{j}, \leq_{j}, \tau_{j})$, $j\in \N\cup\{\infty\}$, is 
 isomorphically embedded in it, i.e. there exist inclusion  maps $\iota_{j}: X_{j}\hookrightarrow X$  such that  for every $x^{1}_{j}, x^{2}_{j}\in X_{j}$,  for every $j\in \N\cup\{\infty\}$, the following holds:
\begin{itemize}
\item ${\sfd}  (\iota_{j}(x^{1}_{j}), \iota_{j}(x^{2}_{j}))= \sfd_{j} (x^{1}_{j}, x^{2}_{j})$;
\item  $x^{1}_{j} \leq_{j} x^{2}_{j}$ 
if and only if $\iota_{j}(x^{1}_{j}) \leq  \iota_{j}(x^{2}_{j})$; 
\item $\tau  (\iota_{j}(x^{1}_{j}), \iota_{j}(x^{2}_{j}))= \tau_{j} (x^{1}_{j}, x^{2}_{j})$.
\end{itemize}
\item The measures $(\iota_{j})_{\sharp} \mm_{j}$  converge to $(\iota_{\infty})_{\sharp} \mm_{\infty}$ weakly in duality with $C_{c}( X)$ in $ X$, i.e.
\begin{equation}\label{eq:weakconv}
\int \varphi\; (\iota_{j})_{\sharp} \mm_{j} \to \int \varphi \; (\iota_{\infty})_{\sharp} \mm_{\infty}, \quad \text{for all }\varphi \in C_{c}( X),
\end{equation} 
where  $C_{c}( X)$ denotes the set of continuous functions with compact support.
\item Convergence of reference points: $\iota_j(\star_j)\to \iota_\infty(\star_\infty)$ in $X$.
\item Volume non-collapsing: there exists a function $v:(0,\infty)\to (0,\infty)$ such that for every $x_{j}\in X_{j}$ it holds $\mm_{j}(B^{\sfd_{j}}_{r}(x_{j})) \geq v(r)>0$.
\item There exists a function $\omega:[0,\infty)\to [0,\infty)$ with $\lim_{r\downarrow 0} \omega(r)=0$ and there exist $p\in (0,1), K\in \R$ such that  $(X_{j},\sfd_{j}, \mm_{j}, \ll_{j}, \leq_{j}, \tau_{j})$ has time-like Ricci curvature bounded above by $K$ with respect to $p\in (0,1)$, $r_{0}>0$ and with remainder function $\omega$ in the synthetic sense of Definition \ref{def:RUB}.
\end{enumerate}
Then also the limit space 
$(X_{\infty},\sfd_{\infty}, \mm_{\infty}, \ll_{\infty}, \leq_{\infty}, \tau_{\infty})$  
 has time-like Ricci curvature bounded above by $K$ with respect to $p\in (0,1), r_{0}+1$ and with remainder function $\omega$ in the synthetic sense of Definition \ref{def:RUB}.
\end{theorem}

The proof can be performed in the same spirit as the proof of Theorem \ref{Thm:WeakStabTCD}, the interested reader is referred to \cite[Thm B.6]{MoSu}.

 \subsection{Synthetic Einstein's vacuum equations}\label{SS:SE}
 Combining the synthetic upper and lower bounds on the time-like Ricci curvature, i.e. Definitions \ref{def:TCD(KN)} and \ref{def:RUB},  the following synthetic version for the  vacuum Einstein equations (with possibly non-zero cosmological constant) is rather natural.
 
  \begin{definition}[Synthetic vacuum Einstein's equations]\label{def:SVEE}
  Fix $p\in (0,1)$, $\Lambda\in \R$, $N\in (0,\infty]$.\ We say that the measured  Lorentzian  pre-length space $(X,\sfd,\mm, \ll, \leq, \tau)$ satisfies the (resp. weak) synthetic formulation of the vacuum Einstein equations $\Ric\equiv \Lambda$ with cosmological constant $\Lambda\in \R$ and has synthetic dimension $\leq N$ if 
\begin{itemize}
  \item  $(X,\sfd,\mm, \ll, \leq, \tau)$ satisfies the $\TCD^{e}_{p}(\Lambda,N)$ (resp. $\wTCD^{e}_{p}(\Lambda,N)$) condition;
  \item There exists $r_{0}>0$ and a function  $\omega:[0,r_{0})\to [0,\infty)$ with $\lim_{r\downarrow 0} \omega(r)=0$ such that $(X,\sfd,\mm, \ll, \leq, \tau)$ has time-like Ricci curvature bounded above by $\Lambda$, with respect to  $p\in (0,1), r_{0}$ and $\omega$.
\end{itemize} 
\end{definition}
 
 Combining the stability of time-like Ricci lower and upper bounds (i.e. Theorem \ref{Thm:WeakStabTCD} and Theorem \ref{thm:stabUB}) gives the stability of the synthetic vacuum Einstein equations under the aforementioned natural Lorentzian variant of measured Gromov-Hausdorff convergence. 
 
\begin{theorem}[Weak stability of synthetic vacuum Einstein's equations]\label{thm:StabVEE}
Let  $\{(X_{j},\sfd_{j}, \mm_{j}, \star_{j}, \ll_{j}, \leq_{j}, \tau_{j})\}_{j\in \N\cup\{\infty\}}$ be a sequence of pointed measured Lorentzian geodesic spaces  satisfying the following properties:
\begin{enumerate}
\item There exists a locally causally closed, $\cK$-globally hyperbolic  Lorentzian geodesic space $(X,  \sfd,   \ll,  \leq,  \tau)$ such that each $(X_{j},\sfd_{j}, \mm_{j}, \ll_{j}, \leq_{j}, \tau_{j})$, $j\in \N\cup\{\infty\}$, is isomorphically embedded in it (as in ${\rm (1)}$ of Theorem \ref{thm:stabUB}).
\item The measures $(\iota_{j})_{\sharp} \mm_{j}$  converge to $(\iota_{\infty})_{\sharp} \mm_{\infty}$ weakly in duality with $C_{c}( X)$ in $ X$, i.e. \eqref{eq:weakconv} holds.
\item Convergence of reference points: $\iota_j(\star_j)\to \iota_\infty(\star_\infty)$ in $X$.
\item  Volume non-collapsing: there exists a function $v:(0,\infty)\to (0,\infty)$ such that for every $x_{j}\in X_{j}$ it holds $\mm_{j}(B^{\sfd_{j}}_{r}(x_{j})) \geq v(r)>0$.\
\item There exist $p\in (0,1), \Lambda \in \R, N\in (0,\infty], r_{0}>0$ and $\omega:[0,r_{0})\to [0,\infty)$ with $\lim_{r\downarrow 0} \omega(r)=0$ such that   $(X_{j},\sfd_{j}, \mm_{j}, \ll_{j}, \leq_{j}, \tau_{j})$  satisfies the synthetic formulation of the vacuum Einstein equations $\Ric\equiv \Lambda$ with cosmological constant $\Lambda\in \R$, with synthetic dimension $\leq N$ with respect to $p\in (0,1), r_{0}$ and $\omega$ as in Definition \ref{def:SVEE}.
\end{enumerate}
Then  the limit space $(X_{\infty},\sfd_{\infty}, \mm_{\infty}, \ll_{\infty}, \leq_{\infty}, \tau_{\infty})$  satisfies the weak synthetic formulation of the vacuum Einstein equations $\Ric\equiv \Lambda$ with cosmological constant $\Lambda\in \R$, with synthetic dimension $\leq N$ with respect to $p\in (0,1), r_{0}+1$ and $\omega$ as in Definition \ref{def:SVEE}.
\end{theorem}

\begin{remark}
By \cite[Thm 4.7, Rem. 4.8]{MoSu} and  Theorem \ref{T:smoothTCD}, if $(X,\sfd,\mm, \ll, \leq, \tau)$ is a (for simplicity say a compact subset in a) smooth Lorentzian manifold, then $(X,\sfd,\mm, \ll, \leq, \tau)$ satisfies the Einstein equations $\Ric\equiv \Lambda$ in the smooth classical sense if and only if $(X,\sfd,\mm, \ll, \leq, \tau)$ satisfies the Einstein equations in the synthetic sense of Definition   \ref{def:SVEE}.\  
Therefore, Theorem \ref{thm:StabVEE} gives that the corresponding limits of smooth solutions to Einstein's equation  $\Ric\equiv \Lambda$  satisfy the weak synthetic Einstein's equations $\Ric\equiv \Lambda$ in the sense of Definition   \ref{def:SVEE}.\  In other words, the vacuum Einstein equations are stable under the conditions (and with respect to the notion of convergence) of   Theorem \ref{thm:StabVEE}.
\smallskip

Let us mention that the stability of the Einstein equations under various notions of (weak) convergence is a subject of high interest in general relativity.
\\Classically, the problem is stated in terms of convergence of a sequence of Lorentzian metrics $g_j$ converging to a limit Lorentzian metric $g_{\infty}$, \emph{on a fixed underlying manifold}.
\\If $g_j$  solve the vacuum Einstein equations,  $g_j\to g_\infty$ in $C^{0}_{loc}$ and the derivatives of $g_j$ converge in $L^2_{loc}$, then it is well known that the limit $g_{\infty}$ satisfies the vacuum Einstein equations as well.
\\However, if the $g_j\to g_\infty$ in $C^{0}_{loc}$ and the derivatives of $g_j$ converge \emph{weakly} in $L^2_{loc}$, explicit examples have been constructed  (see \cite{Burn,GrWa} for examples  in symmetry classes) where the limit $g_{\infty}$ may satisfy the Einstein equations with a \emph{non-vanishing} stress energy momentum tensor.\  Burnett \cite{Burn} conjectured that, if there exist $C>0$ and $\lambda_{j}\to 0$ such that
\[
|g_{j}- g_{\infty}| \leq \lambda_{j}, \quad |\partial g_{j}|\leq C, \quad |\partial^{2} g_{j}| \leq C \lambda_{j}^{-1},
\]
then $g_{\infty}$ is isometric to a solution to the Einstein-massless Vlasov system for some suitable Vlasov field.\ Such a conjecture remains open,  although there has been recent progress \cite{HuLuk1, HuLuk2,GTe} under symmetry conditions.\ We also mention the recent work \cite{LukRod} where concentrations (at the level of $\partial g_{j}$) are allowed in addition to oscillations.
\smallskip

Theorem \ref{thm:StabVEE} gives a new point of view on the stability of the vacuum Einstein equations.\ Indeed, while in the aforementioned results the  metrics $g_{j}$ are converging \emph{on a fixed underlying manifold}, in Theorem \ref{thm:StabVEE} also the underlying space $X$ may vary (along the sequence and in the limit), allowing change in topology in the limit, as one may expect in case of formation of singularities.\ Moreover, the notion of convergence is quite different in spirit: while in the aforementioned results $g_{j}\to g_{\infty}$ in a suitable \emph{functional analytic sense}, in Theorem \ref{thm:StabVEE} the spaces are converging in a \emph{more geometric sense}  (inspired by the pointed measured Gromov-Haudorff convergence).
\end{remark}


\section{Some possible research directions and open problems}\label{Sec:Open}
Due to the richness of the theory of metric measure spaces satisfying the Lott-Sturm-Villani  Curvature-Dimension condition $\CD(K,N)$, it is natural to expect that the Lorentzian counterpart surveyed here will enjoy a rich theory as well.\ In this section we propose some possible research directions and open problems, without any attempt to be exhaustive; on the contrary, we believe that Poincar\'e's famous sentence ``surprising results shall be obtained'' is as appropriate as ever.

\subsection*{Convergence of spaces and pre-compactess}
   A fundamental property of the class of $\CD(K,N)$ spaces (say with a bound on the diameter for simplicity) is the compactness under measured Gromov-Hausdorff (mGH for short) convergence.\ This is a consequence of two deep results: the Gromov pre-compacteness theorem \cite{Gro} and the stability of the $\CD(K,N)$ condition under mGH convergence \cite{sturm:I, sturm:II, lottvillani:metric}.
    
     In Sections \ref{sec:stabLB}, \ref{SS:UB}, \ref{SS:SE}, we proposed  possible notions of convergence for Lorentzian synthetic spaces and proved stability results for synthetic Ricci bounds (and for a synthetic characterization of the Einstein vacuum equations), thus ticking the latter key property.
    
     It is an open problem whether the class of $\TCD^e_p(K,N)$ spaces (say with bounded diameter) satisfy some (pre)-compactness result in the spirit of the Gromov pre-compactness theorem (with respect to the convergence proposed in Definition \ref{D:convergence}, or some other suitable notion of convergence). 
     
     A key difficulty is that while the $\CD/\MCP$ conditions for metric measure spaces imply a control on the volume growth of metric balls (more precisely a local volume doubling property), in our setting the $\tau$-balls typically have infinite volume (e.g., in Minkowski space, $\tau$-spheres are  hyperboloids).\ Thus  we cannot expect to have compactness in classical pointed-measured-Gromov-Hausdorff topology (which is thus not anymore the clearly natural notion for weak convergence of spaces).

\subsection*{Eulerian vs Lagrangian approach}
The point of view on synthetic Ricci bounds employed in this presentation has been \emph{Lagrangian}, i.e. we analysed the convexity/concavity properties of the entropy along suitable Lorentz-Wasserstein geodesics.\ In the smooth setting (as well in the metric measure setting, though the non-smoooth framework is more delicate) such a Lagrangian point of view is equivalent to an \emph{Eulerian} approach based on the Bochner inequality (which is also known as Bakry-\'Emery $\Gamma$-calculus), see \cite{AGS1, AGS2, AGMR, EKS,AMS}.\ It would be very interesting to develop an Eulerian approach to timelike Ricci curvature bounds, and possibly prove the equivalence with the Lagrangian point of view (under suitable assumptions, for instance a natural Lorentzian counterpart of the assumption that the Cheeger energy is a quadratic form).

Building such a bridge between Eulerian and Lagrangian approaches will very likely open the door to a Lorentzian theory for gradient flows, which has been very rich in the smooth Riemannian (and non-smooth metric-measure) setting (see for instance the monograph \cite{AGSBook}, the aforementioned \cite{AGS1, AGS2, AGMR, EKS,AMS} and the more recent survey \cite{AmbICM}).

\subsection*{Compatibility with synthetic sectional curvature lower bounds}
In the smooth setting, it is clear that non-negative (or more generally a lower bound on the) sectional curvature  implies non-negative (or more generally a lower bound on the) Ricci curvature, as the latter is obtained as a trace of the former. 

Synthetic lower bounds on the timelike sectional curvature in terms of comparison triangles (\'a la Alexandrov) for Lorentzian synthetic spaces have been proposed and studied in \cite{KS} (see also the more recent \cite{BeSa}).\ It is a natural open problem whether such synthetic lower bounds on the timelike sectional curvature imply the synthetic lower bounds on the timelike Ricci curvature surveyed here. Note that while the former is purely metric, the latter needs to fix a  reference volume measure in order to be formulated.\ Natural candidates are the canonical volume measures on Lorentzian (pre-)length spaces constructed in \cite{McSa}.

\subsection*{Null synthetic Ricci curvature lower bounds}
In the present survey, we only treated \emph{timelike} Ricci curvature lower bounds.\ This was motivated by the applications we had in mind, namely Hawking's singularity theorem, time-like Bishop-Gromov, time-like Bonnet Myers (see Section \ref{S:applications}, or \cite{CaMo:20} for a more complete list) as well as a synthetic formulation of Einstein's equations (see Section \ref{S:SEE}, or \cite{MoSu} for an optimal transport characterization of the full Einstein equations possibly with a non-zero stress-energy tensor).

It is however also an interesting direction to find a synthetic characterization of  \emph{null} Ricci curvature lower bounds, i.e. $\Ric_g(v,v)\geq K$ for all \emph{null} vectors $v\in TM$.\ Indeed this is at the heart of Penrose's work in general relativity: it is one of the assumptions of his celebrated singularity theorem \cite{Pen} which was awarded the Nobel prize, as well as in much of his work (see for instance \cite[Sec. 1.9, 6.3, 6.4]{La} or \cite{La2}).

\subsection*{Independence of the theory on $p\in (0,1]$}
Note that the definition of the $\TCD^e_p(K,N)$ condition (Definition \ref{def:TCD(KN)}) depends on the exponent $p\in (0,1]$ chosen to fix the optimal transport cost and hence the geodesic structure on the space of probability measures.\ In the metric measure setting (for essentially non-branching metric measure spaces) it was recently proved in \cite{ACCMS} that the $\CD(K,N)$ theory is independent of the  exponent chosen to metrize the space of probability measures.\  Arguing along the lines of \cite{ACCMS}, which in turns builds on \cite{biacava:streconv, cava:MongeRCD, CMi, CM1}, it should be possible to prove the analogous independence of $p$ in the Lorentzian synthethic setting (under suitable timelike non-branching and causality assumptions).

\subsubsection*{Data availability}
The topics of the survey are theoretical and there are no particular data related.\ Detailed bibliographical references have been provided throughout the survey when appropriate. 


%
\end{document}